\documentclass[12pt]{amsart}

\usepackage[latin1]{inputenc}
\usepackage{amsbsy}
\usepackage{amscd} 
\usepackage{amsfonts}
\usepackage{amsgen} 
\usepackage{amsmath}
\usepackage{amsopn} 
\usepackage{amssymb}
\usepackage{amstext}
\usepackage{amsthm} 
\usepackage{amsxtra}
\usepackage[all]{xy}
\usepackage[mathscr]{eucal}
\usepackage{rotating}
\usepackage{MnSymbol}

\theoremstyle{plain} 
\newtheorem{thm}{Theorem}[section]
\newtheorem{prop}[thm]{Proposition}
\newtheorem{lem}[thm]{Lemma}
\newtheorem{cor}[thm]{Corollary}

\theoremstyle{definition}
\newtheorem{defn}[thm]{Definition}
\newtheorem{rem}[thm]{Remark}

\newtheorem{concl}[thm]{Conclusion}

\numberwithin{equation}{section}

\renewcommand{\theta}{\vartheta}
\renewcommand{\phi}{\varphi}
\renewcommand{\epsilon}{\varepsilon}
\renewcommand{\subset}{\subseteq}
\renewcommand{\supset}{\supseteq}

\newcommand{\qqed}{\qed\medskip}
\newcommand{\N}{\mathbb N}
\newcommand{\Z}{\mathbb Z}

\newcommand{\R}{\mathbb R}
\newcommand{\C}{\mathbb C}
\newcommand{\F}{\mathbb F}

\newcommand{\univ}[2]{C^*\left(#1 \; \big| \; #2\right)}

\newcommand{\Proof}{\textbf{Proof. }}

\newcommand{\surjlinks}{\twoheadleftarrow}

\newcommand{\Sn}{S_n^+}
\newcommand{\Snp}{{S_n'}^{\!+}}
\newcommand{\Bn}{B_n^+}
\newcommand{\Bnp}{{B_n'}^{\!+}}
\newcommand{\Bns}{B_n^{\#+}}
\newcommand{\Hn}{H_n^+}
\newcommand{\On}{O_n^+}
\newcommand{\Bnhl}{B_n^{\#*}}
\newcommand{\Hnhl}{H_n^*}
\newcommand{\Onhl}{O_n^*}
\newcommand{\Hnsrund}{H_n^{(s)}}
\newcommand{\Hnseckig}{H_n^{[s]}}
\newcommand{\ASn}{A_{\Sn}}
\newcommand{\ASnp}{A_{\Snp}}
\newcommand{\ABn}{A_{\Bn}}
\newcommand{\ABnp}{A_{\Bnp}}
\newcommand{\ABns}{A_{\Bns}}
\newcommand{\AHn}{A_{\Hn}}
\newcommand{\AOn}{A_{\On}}
\newcommand{\ABnhl}{A_{\Bnhl}}
\newcommand{\AHnhl}{A_{\Hnhl}}
\newcommand{\AOnhl}{A_{\Onhl}}

\newcommand{\idpart}{|}
\newcommand{\paarpart}{\sqcap}
\newcommand{\baarpartbaustein}{\rotatebox{180}{$\sqcap$}}
\newcommand{\baarpart}{
\mathrel{\vcenter{\offinterlineskip \hbox{$\baarpartbaustein$}}}}
\newcommand{\upsubset}{\begin{rotate}{90}$\subset$\end{rotate}}
\newcommand{\downsubset}{\begin{turn}{270}$\subset$\end{turn}}
\newcommand{\singleton}{\uparrow}
\newcommand{\doublesingletonrot}{
\mathrel{\vcenter{\offinterlineskip
\hbox{$\shortmid$} \vskip +.7ex \hbox{$\shortmid$}}}}

\newcommand{\vierpart}{
\mathrel{\offinterlineskip
\hskip0ex\hbox{$\sqcap$}\hskip -.4ex\hbox{$\sqcap$} \hskip -0.4ex\hbox{$\sqcap$}}}

\newcommand{\vierpartrot}{
\mathrel{\vcenter{\offinterlineskip
\hbox{$\baarpart$} \vskip -.1ex \hbox{$\shortmid$} \vskip -.1ex \hbox{$\paarpart$}}}}

\newcommand{\dreipart}{
\mathrel{\offinterlineskip
\hskip0ex\hbox{$\sqcap$}\hskip -.4ex\hbox{$\sqcap$}}}

\newcommand{\crosspart}{
\mathrel{\offinterlineskip
\hbox{$/$}\hskip -.95ex\hbox{$\backslash$}}}

\newcommand{\fatcrosspart}{
\mathrel{\offinterlineskip
\hbox{$\crosspart$}\hskip -.75ex\hbox{$\crosspart$}}}

\newcommand{\midmid}{
\mathrel{\vcenter{\offinterlineskip
\hbox{$\shortmid$} \vskip -1.6ex \hbox{$\shortmid$}}}}

\newcommand{\halflibpart}{
\mathrel{\offinterlineskip
\hbox{$\bigtimes$}\hskip -1.55ex\hbox{$\midmid$}}}

\newcommand{\liegeblkn}{\begin{turn}{270}$[$\end{turn}}

\newcommand{\liegebalken}{
\mathrel{\vcenter{\offinterlineskip
\vskip -1.3ex \hbox{$\liegeblkn$}}}}

\newcommand{\longpr}{
\mathrel{\offinterlineskip
\hskip0ex\hbox{$\shortmid$}\hskip -.8ex\hbox{$\liegebalken$} \hskip -0.75ex\hbox{$\shortmid$}}}

\newcommand{\hochlongpr}{
\mathrel{\vcenter{\offinterlineskip
\vskip -.5ex \hbox{$\longpr$}}}}

\newcommand{\longpair}{
\mathrel{\offinterlineskip
\hskip0ex\hbox{$\shortmid$}\hskip -1.4ex\hbox{$\hochlongpr$} \hskip -1.4ex\hbox{$\shortmid$}}}

\newcommand{\legpart}{
\mathrel{\offinterlineskip
\hskip0ex\hbox{$\shortmid$}\hskip +.8ex\hbox{$\shortmid$} \hskip -2.5ex\hbox{$\longpair$}}}
\newcommand{\hochshortmid}{
\mathrel{\vcenter{\offinterlineskip
\vskip -1.2ex \hbox{$\shortmid$}}}}

\newcommand{\tiefshortmid}{
\mathrel{\vcenter{\offinterlineskip
\vskip +1.2ex \hbox{$\shortmid$}}}}

\newcommand{\legpartrot}{
\mathrel{\offinterlineskip
\hbox{$\tiefshortmid$}\hbox{$\idpart$}\hbox{$\hochshortmid$}}}

\usepackage{calc}
\newcounter{PartitionDepth}
\newcounter{PartitionLength}
\newcommand{\parti}[2]{
 \begin{picture}(#2,#1)
 \setcounter{PartitionDepth}{-1-#1}
 \put(#2,\thePartitionDepth){\line(0,1){#1}}
 \end{picture}}
\newcommand{\partii}[3]{
 \begin{picture}(#3,#1)
 \setcounter{PartitionLength}{#3-#2}
 \setcounter{PartitionDepth}{-1-#1}
 \put(#2,\thePartitionDepth){\line(0,1){#1}}     
 \put(#3,\thePartitionDepth){\line(0,1){#1}}
 \put(#2,\thePartitionDepth){\line(1,0){\thePartitionLength}}
 \end{picture}}

\newcommand{\upparti}[2]{
 \begin{picture}(#2,#1)
 \setcounter{PartitionDepth}{#1}
 \put(#2,0){\line(0,1){#1}}
 \end{picture}}
\newcommand{\uppartii}[3]{
 \begin{picture}(#3,#1)
 \setcounter{PartitionLength}{#3-#2}
 \setcounter{PartitionDepth}{#1}
 \put(#2,0){\line(0,1){#1}}     
 \put(#3,0){\line(0,1){#1}}
 \put(#2,\thePartitionDepth){\line(1,0){\thePartitionLength}}
 \end{picture}}
\newcommand{\uppartiii}[4]{
 \begin{picture}(#4,#1)
 \setcounter{PartitionLength}{#4-#2}
 \setcounter{PartitionDepth}{#1}
 \put(#2,0){\line(0,1){#1}}
 \put(#3,0){\line(0,1){#1}}
 \put(#4,0){\line(0,1){#1}}
 \put(#2,\thePartitionDepth){\line(1,0){\thePartitionLength}} 
 \end{picture}}
\newcommand{\uppartiv}[5]{
 \begin{picture}(#5,#1)
 \setcounter{PartitionLength}{#5-#2}
 \setcounter{PartitionDepth}{#1}
 \put(#2,0){\line(0,1){#1}}
 \put(#3,0){\line(0,1){#1}}
 \put(#4,0){\line(0,1){#1}}
 \put(#5,0){\line(0,1){#1}}
 \put(#2,\thePartitionDepth){\line(1,0){\thePartitionLength}} 
 \end{picture}}

\newcommand{\uppartviii}[9]{
 \begin{picture}(#9,#1)
 \setcounter{PartitionLength}{#9-#2}
 \setcounter{PartitionDepth}{#1}
 \put(#2,0){\line(0,1){#1}}
 \put(#3,0){\line(0,1){#1}}
 \put(#4,0){\line(0,1){#1}}
 \put(#5,0){\line(0,1){#1}}
 \put(#6,0){\line(0,1){#1}}
 \put(#7,0){\line(0,1){#1}}
 \put(#8,0){\line(0,1){#1}}
 \put(#9,0){\line(0,1){#1}}
 \put(#2,\thePartitionDepth){\line(1,0){\thePartitionLength}} 
 \end{picture}}
\begin{document}
\title{On the classification of easy quantum groups}
\author{Moritz Weber}
\address{Saarland University, Fachbereich Mathematik, Postfach 151150,
66041 Saarbr\"ucken, Germany}
\email{weber@math.uni-sb.de}
\date{\today}
\maketitle

\begin{abstract}
In 2009, Banica and Speicher began to study the compact quantum groups $G$ with $S_n\subset G\subset \On$ whose intertwiner spaces are induced by some partitions. These so-called \emph{easy} quantum groups have a deep connection to combinatorics. We continue their work on classifying these objects, by introducing some new examples of easy quantum groups. In particular, we show that the six easy groups $O_n$, $S_n$, $H_n$, $B_n$, $S_n'$ and $B_n'$ split into seven cases $\On$, $\Sn$, $\Hn$, $\Bn$, $\Snp$, $\Bnp$ and $\Bns$ on the side of free easy quantum groups. Also, we give a complete classification in the half-liberated and in the non-hyperoctahedral case.
\end{abstract}

\section*{Introduction}

In 2009, Banica and Speicher (\cite{BS09}) initiated the classification program for easy quantum groups. These are compact quantum groups $G$ such that $S_n\subset G\subset \On$, where $S_n$ is the symmetric group and $\On$ is the (free) orthogonal quantum group, constructed by Wang \cite{W1}. If, in addition, the intertwiner space of $G$ (which determines the quantum group $G$ by Woronowicz, \cite{Wo}) is induced by certain partitions, the quantum group $G$ is called \emph{easy}. In other words, in order to understand the quantum groups in between the symmetric quantum group $\Sn$ (by Wang \cite{W2}) and $\On$, or -- more generally -- in between $S_n$ and $\On$, Banica and Speicher restrict their attention to those that are given by some underlying combinatorial data. 

Furthermore, the class of easy quantum groups appears very natural, since the intertwiner spaces of the (classical) groups $S_n$ and $O_n$ are given by all partitions resp. by all pair partitions, whereas their free analogues $\Sn$ and $\On$ correspond to all noncrossing partitions resp. to all noncrossing pair partitions. Thus, the connection of the quantum groups in between $S_n$ and $\On$ to combinatorics is not at all artificial, quite the contrary, it is an intrinsic background that needs to be exploited. 
Besides, it is a well known feature in free probability theory that the step from ``commutative'' to ``free'' is reflected by passing from all partitions to the noncrossing ones. So, also from this perspective it is worth to study this effect in the present case.

The easy quantum groups have been investigated by Banica, Speicher, Curran, Bichon, Collins and some others in a couple of papers (see for instance \cite{BS09}, \cite{BCS10}, \cite{BCS11}, \cite{BBC07}). Further operator algebraic aspects of some of the easy quantum groups have been studied in \cite{VV07}, \cite{B12a}, \cite{B12b}, \cite{F12}.
The link between easy quantum groups and free probability theory has been discovered by K\"ostler and Speicher, \cite{KS09}. Their noncommutative version of the De Finetti theorem hints at the deep connection between co-actions of (easy) quantum groups and free probability theory, which lead to further work by Banica, Curran and Speicher (see also \cite{C09}, \cite{C10}, \cite{C11}, \cite{CS11}, \cite{BCS12}).

In the case of classical easy quantum groups (which are in fact groups) \linebreak $S_n\subset G\subset O_n$,  the classification is complete and has been done in \cite{BS09}. There are exactly six easy groups: $O_n$, $H_n$, $S_n$, $B_n$, $S_n'$ and $B_n'$. In the case of free easy quantum groups $\Sn\subset G\subset \On$, one quantum group was missing on the list in \cite{BS09}. We very much follow the proof of Banica and  Speicher and we show that there are exactly seven free easy quantum groups (also called \emph{free orthogonal quantum groups}): $\On$, $\Hn$, $\Sn$, $\Bn$, $\Snp$, $\Bnp$ and $\Bns$. The subtlety in the free case is that the group $B_n'$ splits into two cases $\Bnp$ and $\Bns$, whereas all other easy groups and free easy quantum groups are in a natural one-to-one correspondence. We describe this phenomenon to a certain extent although it is still not completely understood. Nevertheless, this shows again that the ``free'' world is somewhat richer than the commutative one.

In \cite{BS09} and \cite{BCS10}, further examples of easy quantum groups are given, the half-liberated versions $\Onhl$ and $\Hnhl$ of $\On$ resp. of $\Hn$ as well as the (infinite) hyperoctahedral series $\Hnsrund$ and $\Hnseckig$. We extend the list by a half-liberated version $\Bnhl$ of $\Bns$ and a further one. We also prove that there are exactly three half-liberated easy quantum groups besides the hyperoctahedral series $\Hnsrund$. Hence, the classification is complete in the classical, in the free and in the half-liberated case. However, the complete classification of all easy quantum groups still remains open.

The paper is organized as follows. In Section \ref{Sect1} we give a short introduction into the concept of the classification of easy quantum groups. For details on this, we refer to the initial paper by Banica and Speicher \cite{BS09}. In Section \ref{Sect2} we focus on the classification of the free easy quantum groups $\Sn\subset G\subset \On$ following the path of \cite{BS09}. We  formulate it in a language of generating partitions, which might help in the future work on the classification. The main results of this section are the description of all possible categories of noncrossing partitions 
and the proof that there are exactly seven of them,
which has partly been done by \cite{BS09}.

A conceptual explanation for the difference of the two quantum versions of $B_n'$ is given in Section \ref{Sect3}, where we study the $C^*$-algebraic level of these objects. We show that the $C^*$-algebras associated to $\Bnp$ and $\Bns$ may be constructed out of the one of $\Bn$ by a tensor product resp. a free product with $C^*(\Z_2)$. Likewise, $\Snp$ is given by a tensor product of $\Sn$ with $C^*(\Z_2)$. We  use a result on the $K$-theory of $\On$ by Voigt (\cite{V11}) to compute the $K$-theory of $\Bn$, $\Bnp$ and $\Bns$. 
The computation of the $K$-groups for $\Sn$, $\Snp$ and $\Hn$ has to be left open.
Finally, we prove that none of the full $C^*$-algebras corresponding to the seven free easy quantum groups is exact, if $n\geq 5$.

Section \ref{Sect4} is devoted to the study of the half-liberated versions $\Onhl$, $\Hnhl$ and $\Bnhl$, the latter one being new. We describe their categories and prove that these are the only half-liberated easy quantum groups, together with the series $\Hnsrund$. We recover a result by Banica, Curran and Speicher (\cite{BCS10}) describing all nonhyperoctahedral easy quantum groups and we show that there are exactly 13 of them (including the two new examples of our article).

In Section \ref{Sect5} we compute the laws of the character $\sum u_{ii}$ of the new quantum groups. By this, we give a second explanation for the existence of two quantum versions of $B_n'$. The laws of the quantum groups $\Bnhl$ and  $\Bns$ are ``complex'' versions of the laws of $B_n'$ and $\Bnp$.

In the last section, we introduce an additional example of an easy quantum group. 

\section{The concept of easy quantum groups}\label{SectConcept}\label{Sect1}

The quantum groups $\Sn$ and $\On$
introduced by Wang (\cite{W1}, \cite{W2}) are milestones in the theory of quantum groups. Hence, it is natural to study the quantum groups lying in between, the so-called homogeneous quantum groups. More generally, the quantum groups in between 
the symmetric group $S_n$ and $\On$ are to be investigated.

In this section we give an overview on the idea of the classification of easy quantum groups, which is due to Banica and Speicher and we mainly follow their article from 2009 (\cite{BS09}). Details may be found there.

\begin{defn} \label{DefOrthQG}
 An \emph{orthogonal} quantum group $G$ is given by a $C^*$-algebra $A$ generated by $n^2$ self-adjoint generators $u_{ij}$, $i,j=1,\ldots, n$ fulfilling the relations
\[\sum_{k=1}^n u_{ik}u_{jk}=\sum_{k=1}^n u_{ki}u_{kj}=\delta_{ij} \quad \forall i,j=1,\ldots,n\]
In other words, the $n\times n$-matrix $u$ formed by the elements $u_{ij}$ is orthogonal, i.e. $u^tu=uu^t=1$.

Furthermore, we require the existence of the following homomorphisms:
\begin{itemize}
 \item $\Delta:A\to A\otimes A$, mapping $u_{ij}$ to $\sum_{k=1}^n u_{ik}\otimes u_{kj}$
\item $\epsilon:A\to \C$, mapping $u_{ij}$ to $\delta_{ij}$
\item $S:A\to A^{op}$, mapping $u_{ij}$ to $u_{ji}$
\end{itemize}
By this definition, $A$ endowed with $\Delta,\epsilon$ and $S$ is a compact quantum group, i.e. the maps $\Delta, \epsilon$ and $S$ satisfy the conditions of a comultiplication, a counit and an antipode. We write $A=A_G$.

In other words, a compact quantum group $G$ is orthogonal, if $G\subset \On$.
\end{defn}

Typically, $A_G$ is a universal $C^*$-algebra generated by $n^2$ generators $u_{ij}$ satisfying the orthogonality relations and some further ones. If $G\subset O_n$ is a group, the algebra $A_G=\mathcal C(G)$ of continuous functions over $G$ equipped with $\Delta, \epsilon$ and $S$ is an orthogonal quantum group. Here, the generators $u_{ij}$ of $\mathcal C(G)$ commute. Thus, any subgroup of the orthogonal group may be seen as an orthogonal quantum group, whereas the converse is not true: In general, a quantum group is not a group.

Philosophically, we could see the $C^*$-algebra $A_G$ of the preceding definition as the ``noncommutative continuous functions'' over $G$. In this sense, a ($C^*$-algebraic compact) quantum group is a \emph{virtual} object given by its algebra of functions and a comultiplication, a counit and an antipode. This duality between a ``space'' (with some additional structure) and the algebra of continuous functions over it (with some additional structure) is a very central and general concept in operator algebras and functional analysis (see for instance \cite[chapter 1]{GVF}). Thus, it is natural to consider a quantum group $G$ as a \emph{true} object and to write expressions like ``$S_n\subset G$'', which become practical by the dual formulation $\mathcal C(S_n)\leftarrow A_G$ (the homomorphism respecting the comultiplication etc.) in terms of algebras over $S_n$ and $G$. 

There is a correspondence between orthogonal quantum groups and tensor categories with duals by Woronowicz \cite{Wo}.
Namely, for $n\in \N$ and an orthogonal quantum group $G$ we denote by $C_G$ the collection of its intertwiner spaces:
\[C_G(k,l)=\{T:(\C^n)^{\otimes k}\to(\C^n)^{\otimes l} \textnormal{ linear map} \; | \; Tu^{\otimes k}=u^{\otimes l}T\} \qquad k,l\in \N_0\]
Here, $u$ is the $n\times n$-matrix imposed by the generators $u_{ij}$ of $G$. Also, we use the convention $(\C^n)^{\otimes 0}=\C$. This correspondence between $G$ and $C_G$ is one-to-one. 

For details on this, see \cite[def. 3.5, def. 1.1 and th. 3.6]{BS09}.

\begin{defn}
A \emph{homogeneous} quantum group is an orthogonal quantum group $G$, such that $\Sn\subset G\subset \On$. This means, there are surjections $\ASn\surjlinks A_{G} \surjlinks \AOn$ which compose to the canonical map from $\AOn$ to $\ASn$. Note that $\AOn$ and $\ASn$ are also denoted by $A_o(n)$ and $A_s(n)$ in the literature.
\end{defn}

For homogeneous quantum groups $G$, we can specify the structure of the intertwiner spaces $C_{G}$ in the following way, using functoriality of the above correspondence \cite[th. 3.8, th. 3.9]{BS09}.
\[\textnormal{span}\{T_p\;|\;p\in NC\}=C_{\Sn}\supset C_{G} \supset C_{\On}=\textnormal{span}\{T_p\;|\;p\in NC_2\}\]

Here, $NC$ denotes the set of all noncrossing partitions (or, more precisely, the collection of all noncrossing partitions $NC(k,l)$, for all $k,l\in\N_0$), whereas $NC_2$ is the set of all noncrossing pairings (cf. \cite[def. 3.7, def. 1.5]{BS09}). The combinatorial object $NC(k,l)$ consists of all possible ways of connecting $k$ upper points with $l$ lower points by lines which are not allowed to cross. Of course, upper points may also be connected to other upper points, as well as lower points with lower points. There may also exist points which are not connected to any other points, the so-called \emph{singletons}.

By $T_p$ we denote the linear map induced by a noncrossing partition $p\in NC(k,l)$, or more generally by a -- possibly crossing -- partition $p\in P(k,l)$, given by:
\[T_p(e_{i_1}\otimes \ldots \otimes e_{i_k}) = \sum_{j_1,\ldots,j_l=1}^n \delta_p(i,j) e_{j_1}\otimes \ldots \otimes e_{j_l}\]
Here, $e_1,\ldots,e_n$ is the canonical basis of $\C^n$, and $\delta_p(i,j)=1$ if and only if the indices $i=(i_1,\ldots,i_k)$ coincide with the indices $j=(j_1,\ldots,j_l)$ when connected by the partition $p$. Otherwise $\delta_p(i,j)=0$.
(cf. \cite[def. 1.7, def. 1.6]{BS09})

For homogeneous quantum groups $G$, $C_{G}$ is a linear subspace of $\textnormal{span}\{T_p\;|\;p\in NC\}$, i.e. $C_G(k,l)$ is a linear subspace of $\textnormal{span}\{T_p\;|\;p\in NC(k,l)\}$ for all $k,l\in \N_0$.
This gives rise to the very natural definition of easy quantum groups by Banica and Speicher \cite[def. 3.10, def. 6.1]{BS09}, restricting the class of homogeneous quantum groups to those whose intertwiner spaces are spanned by subsets of $\{T_p\;|\;p\in NC\}$.

\begin{defn}\label{DefFreeEasy}
 A homogeneous quantum group $G$ is called \emph{free easy} (or \emph{free orthogonal quantum group}), if its intertwiner space is of the form
\[C_{G}=\textnormal{span}\{T_p\;|\; p\in NC_X\}\]
where $NC_2\subset NC_X\subset NC$. More precisely, 
\[C_{G}(k,l)=\textnormal{span}\{T_p\;|\; p\in NC_X(k,l)\}\]
where $NC_2(k,l)\subset NC_X(k,l)\subset NC(k,l)$ for all $k,l\in \N_0$.
\end{defn}

Since $C_{G}$ is a tensor category, it can be shown that $NC_X$ reflects exactly this structure (cf. \cite[def. 3.11, prop. 3.12, def. 1.8]{BS09}). In fact, it is a category of noncrossing partitions, which is a special case of a full category of partitions (see \cite[def. 6.3]{BS09}).
In the sequel, $P(k,l)$ denotes the set of possibly crossing partitions connecting $k$ upper points with $l$ lower points.

\begin{defn}\label{DefFullKateg}
 A \emph{(full) category of partitions} $P_X$ is a collection of subsets\linebreak $P_X(k,l)\subset P(k,l)$, such that:
\begin{itemize}
 \item $P_X$ is stable by \emph{tensor product}. 

(If $p\in P(k,l)$ and $q\in P(k',l')$, then $p\otimes q\in P(k+k',l+l')$ is obtained by horizontal concatenation, i.e. the first $k$ of the $k+k'$ upper points are connected by $p$ to the first $l$ of the $l+l'$ lower points, whereas $q$ connects the remaining $k'$ upper points with the remaining $l'$ lower points.)

\item $P_X$ is stable by \emph{composition}. 

(If $p\in P(k,l)$ and $q\in P(l,m)$, then $pq\in P(k,m)$ is obtained by vertical concatenation: connect $k$ upper points by $p$ to $l$ middle points and continue the lines by $q$ to $m$ lower points. This yields a partition, connecting $k$ upper points with $m$ lower points. The middle points $l$ are removed as well as all the  lines that have no connection to the upper or the lower points.)

\item $P_X$ is stable by \emph{involution}.

(If $p\in P(k,l)$, then $p^*\in P(l,k)$ is obtained by turning $p$ upside down.)

\item $P_X$ is stable by \emph{rotation}.

(Let $p\in P(k,l)$ connecting $k$ upper points with $l$ lower points. Shifting the first left upper point to the left of the lower points -- without changing the connections between the points -- gives rise to a partition in $P(k-1, l+1)$. This procedure may be iterated and also be applied in the converse direction, as well as to the right side.
In particular, for a partition $p\in P(0,l)$, we may rotate the very left point to the very right and vice versa.)

\item $P_X$ contains the \emph{unit partition} $\idpart\in NC(1,1)$ connecting one upper point with one lower point.

\item $P_X$ contains the \emph{pair partition} (also called the \emph{duality partition}) \linebreak $\paarpart\in NC(0,2)$ connecting two lower points.
\end{itemize}
We will call these operations (tensor product, composition, involution and rotation)  the \emph{category operations}.
\end{defn}

\begin{defn} \label{DefKatNC}
 A \emph{category of noncrossing partitions} $NC_X$ is a category of partitions, such that $NC_X\subset NC$, i.e. all partitions in $NC_X$ are noncrossing.
\end{defn}

\begin{rem}
 Stability by rotation may be deduced from the other properties (see the proof of lemma 2.7 in \cite{BS09}). If $P_X$ is a category of partitions, rotation allows us to restrict our attention to partitions in $P(0,k)$ for $k\in \N$. We will often make use of this.
\end{rem}

The category operations enable us to ``connect'' two partitions in a category, or to ``erase'' certain points. By this we mean the following. If, for instance, the three block -- i.e. the partition in $NC(0,3)$, connecting three points in a row -- is in a category $P_X$, then also the four block, since it may be constructed in $P_X$ in the following way. (We compose the tensor product of two three block partitions with $\idpart\otimes\idpart\otimes\paarpart^*\otimes\idpart\otimes\idpart$.)
\setlength{\unitlength}{0.5cm}
\begin{center}
\begin{picture}(11,4)
 \put(-0.1,2.5){\uppartiii{1}{1}{2}{3}}
 \put(-0.1,2.5){\uppartiii{1}{4}{5}{6}}
 \put(-0.1,3){\parti{2}{1}}
 \put(-0.1,3){\parti{2}{2}}
 \put(-0.1,3){\partii{1}{3}{4}}
 \put(-0.1,3){\parti{2}{5}}
 \put(-0.1,3){\parti{2}{6}}
 \put(7,1.8){$=$}
 \put(7,1.5){\uppartiv{1}{1}{2}{3}{4}}
 \end{picture}
\end{center}

We have ``connected'' two copies of the three block.
As a second example, let the following partition $p$ be in $P_X$.

\setlength{\unitlength}{0.5cm}
\begin{center}
\begin{picture}(8,4)
 \put(1,0){\uppartiii{3}{1}{6}{7}}
 \put(1,0){\uppartii{2}{2}{5}}
 \put(1,0){\uppartii{1}{3}{4}}
 \put(0,1){$p \;=$}
 \end{picture}
\end{center}

Then, the three block partition is in $P_X$ too, due to the following construction.
 
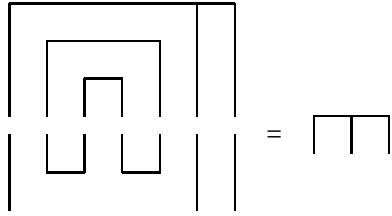
\setlength{\unitlength}{0.5cm}
\begin{center}
\begin{picture}(10,6)
 \put(-0.1,2.5){\uppartiii{3}{1}{6}{7}}
 \put(-0.1,2.5){\uppartii{2}{2}{5}}
 \put(-0.1,2.5){\uppartii{1}{3}{4}}
 \put(-0.1,3){\parti{2}{1}}
 \put(-0.1,3){\partii{1}{2}{3}}
 \put(-0.1,3){\partii{1}{4}{5}}
 \put(-0.1,3){\parti{2}{6}}
 \put(-0.1,3){\parti{2}{7}}
 \put(8,1.8){$=$}
 \put(8,1.5){\uppartiii{1}{1}{2}{3}}
 \end{picture}
\end{center}

We have ``erased'' the points not belonging to the three block in $p$.

The categories of noncrossing partitions arise as intertwiner spaces of compact quantum groups in between $\Sn$ and $\On$, whereas the full categories of partitions correspond to compact quantum groups in between $S_n$ and $\On$. Recall that the permutation group $S_n$ gives rise to the orthogonal quantum group $\mathcal C(S_n)$ in the sense of Definition \ref{DefOrthQG}. Here, the generators consist of the coordinate functions $u_{ij}:S_n\to\C$, mapping $g\in S_n$ to $1$, if the permutation $g$ maps $i\in \{1,\ldots,n\}$ to $j$, and $0$ otherwise.

The homogeneous quantum groups are included in a more general class. It is given by 
all compact quantum groups $G$ with $S_n\subset G\subset\On$. In the dual picture, this means that there are surjective homomorphism $\mathcal C(S_n)\surjlinks A_{G} \surjlinks A_{\On}$ composing to the canonical map from $A_{\On}$ to $\mathcal C(S_n)$.
The consequence for the intertwiner spaces is the following, again by Woronowicz \cite{Wo} (cf. also \cite[th. 1.10]{BS09}):
\[\textnormal{span}\{T_p\;|\;p\in P\}=C_{\mathcal C(S_n)}\supset C_{G} \supset C_{\On}=\textnormal{span}\{T_p\;|\;p\in NC_2\}\]

Here, $P$ is the collection of all (possibly crossing) partitions $P(k,l)$ (see \cite[def. 1.5]{BS09}).
In analogy to Definition \ref{DefFreeEasy}, we can now define the class of easy quantum groups. (cf. \cite[def. 6.1]{BS09}, \cite[def. 2.1]{BCS10})

\begin{defn}\label{DefEasyQG}
 A compact quantum group $G$, with $S_n\subset G\subset \On$ is called \emph{easy}, if its intertwiner space is of the form
\[C_{G}=\textnormal{span}\{T_p\;|\; p\in P_X\}\]
where $NC_2\subset P_X\subset P$. 
\end{defn}

Again, it can be shown that $P_X$ inherits the structure of a tensor category. Thus, $P_X$ is a full category of partitions, see \cite[prop. 6.4]{BS09}.

\begin{concl}
We conclude that the problem of classifying the free easy quantum groups is equivalent to classifying the categories of noncrossing partitions.
\begin{align*}
 \{\textnormal{free easy quantum groups}\} \quad &\longleftrightarrow \quad
\{\textnormal{categories of noncrossing partitions}\}
\end{align*}
More generally, the classification of the easy quantum groups corresponds to the classification of the full categories of partitions. (cf. also \cite[prop. 3.12, prop. 6.4]{BS09})
\begin{align*}
 \{\textnormal{easy quantum groups}\} \quad &\longleftrightarrow \quad
\{\textnormal{full categories of partitions}\}
\end{align*}

\end{concl}

\section{The classification of the free easy quantum groups}\label{Sect2}

For easy quantum groups $\Sn\subset G\subset \On$ a complete classification is feasible. The major part of this has been done by Banica and Speicher in \cite{BS09} even though one quantum group was missing. We fill in this gap using their machinery and we follow the path of their work. We formulate our results in terms of generating partitions which will be explained in the sequel.

First, we define seven orthogonal quantum groups, following partly \cite[def. 3.2, def. 3.4]{BS09}. We will see that these are the only free easy quantum groups.

\begin{defn}\label{DefSieben}
 Let $n\in \N$, let $A$ be a $C^*$-algebra and let $u_{ij}\in A$ be self-adjoint elements, for $i,j=1,\ldots,n$. Denote by $(u_{ij})$ the system of the elements\linebreak $\{u_{ij}, i,j=1,\ldots,n\}$, and by $\univ{u_{ij}, i,j=1,\ldots,n}{\ldots}$ universal $C^*$-algebras.
\begin{itemize}
\item[(1)] $(u_{ij})$ is \emph{orthogonal}, if $\sum_{k=1}^n u_{ik}u_{jk}=\sum_{k=1}^n u_{ki}u_{kj}=\delta_{ij}$ for all  $i,j=1,\ldots,n$.

Let $A_{\On}:= \univ{u_{ij}, i,j=1,\ldots,n}{(u_{ij}) \textnormal{ is orthogonal}}$.

\item[(2)] $(u_{ij})$ is \emph{cubic}, if $(u_{ij})$ is orthogonal and $u_{ik}u_{jk}=u_{ki}u_{kj}=0$ whenever $i\neq j$, for $i,j,k=1,\ldots,n$.

Let $A_{\Hn}:= \univ{u_{ij}, i,j=1,\ldots,n}{(u_{ij}) \textnormal{ is cubic}}$.

\item[(3)] $(u_{ij})$ is \emph{magic'}, if $(u_{ij})$ is cubic and $\sum_{k=1}^n u_{ik}=\sum_{k=1}^n u_{kj}$ for all  $i,j=1,\ldots,n$. We then denote $r:=\sum_{k=1}^n u_{ik}$.

Let $A_{\Snp}:= \univ{u_{ij}, i,j=1,\ldots,n}{(u_{ij}) \textnormal{ is magic'}}$.

\item[(4)] $(u_{ij})$ is \emph{magic}, if $(u_{ij})$ is magic' and $r=1$. 

Let $A_{\Sn}:= \univ{u_{ij}, i,j=1,\ldots,n}{(u_{ij}) \textnormal{ is magic}}$.

\item[(5)]  $(u_{ij})$ is \emph{bistochastic$^\#$}, if $(u_{ij})$ is orthogonal and $\sum_{k=1}^n u_{ik}=\sum_{k=1}^n u_{kj}$ for all \linebreak $i,j=1,\ldots,n$. We then denote $r:=\sum_{k=1}^n u_{ik}$.

Let $A_{\Bns}:= \univ{u_{ij}, i,j=1,\ldots,n}{(u_{ij}) \textnormal{ is bistochastic}^\#}$.

\item[(6)]  $(u_{ij})$ is \emph{bistochastic'}, if $(u_{ij})$ is bistochastic$^\#$ and $u_{ij}r=ru_{ij}$ for all \linebreak $i,j=1,\ldots,n$ and $r=\sum_{k=1}^n u_{ik}$.

Let $A_{\Bnp}:= \univ{u_{ij}, i,j=1,\ldots,n}{(u_{ij}) \textnormal{ is bistochastic'}}$.

\item[(7)]  $(u_{ij})$ is \emph{bistochastic}, if $(u_{ij})$ is bistochastic' and $r=1$.

Let $A_{\Bn}:= \univ{u_{ij}, i,j=1,\ldots,n}{(u_{ij}) \textnormal{ is bistochastic}}$.
\end{itemize}
\end{defn}

Note the change of notation compared to the definition in \cite[def. 3.2, def. 3.4]{BS09}. In our article, the bistochastic' relations include commutativity of $u_{ij}$ and $r$, whereas this is not the case in the more general bistochastic$^\#$ relations. In \cite{BS09}, our bistochastic$^\#$ relations are called bistochastic' and there are no relations, involving the commutativity of the $u_{ij}$ and $r$. In this article, we distinguish between the bistochastic' and the bistochastic$^\#$ relations, which has not been done in \cite{BS09} or in the subsequent papers. For the reasons why we renamed the bistochastic' relations, see also Remark \ref{RemX} and Remark \ref{RemY}.

There are canonical surjections between the above universal $C^*$-algebras, mapping the generators $u_{ij}$ to $u_{ij}$. We illustrate this by the following diagram.
\begin{align*}
 A_{\Bn} &\qquad\longleftarrow &A_{\Bnp} &\qquad\longleftarrow &A_{\Bns} &\qquad\longleftarrow &A_{\On}\\
\downarrow & &\downarrow & & & &\downarrow\\
 A_{\Sn} &\qquad\longleftarrow &A_{\Snp} & &\longleftarrow &  &A_{\Hn}
\end{align*}

\begin{rem}\label{BemUeberR}
 \begin{itemize}
  \item[(a)] If $(u_{ij})$ is bistochastic$^\#$, then $r:=\sum_{k=1}^n u_{ik}$ is a symmetry (i.e. $r=r^*$ and $r^2=1$). Hence, $r$ is also a symmetry in $A_{\Bnp}$ and in $A_{\Snp}$ (and trivially also in $\ABn$ and $\ASn$).
\item[(b)] If $(u_{ij})$ is magic', then $r$ commutes with all $u_{ij}$. Thus, there is no ``magic$^\#$'' version of the magic' relations.
\item[(c)] The system $(u_{ij})$ is magic if and only if $(u_{ij})$ is orthogonal and all $u_{ij}$ are projections, $i,j=1,\ldots,n$.
 \end{itemize}
\end{rem}
\Proof (a) Since the $u_{ij}$ are self-adjoint, so is $r$. Secondly, we write $r^2$ as \linebreak $r^2=\sum_{k=1}^n(u_{ik}r)=\sum_{k=1}^n\sum_{j=1}^n u_{ik}u_{jk}$, using $r=\sum_{j=1}^n u_{jk}$ for all $k$. Now, changing the order of summation and using the orthogonal relations, we end up with \linebreak $r^2 =\sum_{j=1}^n\delta_{ij}=1$.

(b) We compute $u_{ij}r=\sum_{k=1}^n u_{ij}u_{ik}=u_{ij}^2$, using the cubic relations. On the other hand, $ru_{ij}=\sum_{k=1}^n u_{ik}u_{ij}=u_{ij}^2$ by the same argument. By the magic' relations, $r$ can always be written as a sum involving $i$, for all $i=1,\ldots,n$. 

(c) If $(u_{ij})$ is magic, then $u_{ij}^2=\sum_{k=1}^nu_{ij}u_{ik}=u_{ij}\sum_{k=1}^nu_{ik}=u_{ij}$. Conversely, if the $u_{ij}$ are projections and if $(u_{ij})$ is orthogonal, then $\sum_{k=1}^nu_{ik}=\sum_{k=1}^nu_{ik}^2=1=\sum_{k=1}^n u_{kj}$ for all $i,j$. Thus, the projections $u_{ik}$ sum up to a projection, which proves that they are mutually orthogonal.\qqed

\begin{lem}\label{LemSiebenQG}
 The seven $C^*$-algebras of Definition \ref{DefSieben} give rise to orthogonal quantum groups and we write:
\begin{itemize}
 \item $\On$: the free orthogonal quantum group
 \item $\Hn$: the hyperoctahedral quantum group
 \item $\Snp$: the modified symmetric quantum group 
 \item $\Sn$: the symmetric quantum group 
 \item $\Bns$: the freely modified bistochastic quantum group
 \item $\Bnp$: the modified bistochastic quantum group
 \item $\Bn$: the bistochastic quantum group
\end{itemize}
\end{lem}
\Proof Verify that the relations on $(u_{ij})$ are stable by the maps $\Delta$, $\epsilon$ and $S$ of Definition \ref{DefOrthQG}.\qqed

\begin{rem}\label{RemX}
The quantum groups $\On$ and $\Sn$ have been introduced by Wang (\cite{W1}, \cite{W2}), $\Hn$ was defined by Banica, Bichon and Collins in \cite{BBC07}, and $\Snp$ and $\Bn$ were given by Banica and Speicher in \cite{BS09}. For $\Bnp$ and $\Bns$ there is a subtlety to remark. We split the idea of the quantum group $\Bnp$ of Banica and Speicher into two quantum groups $\Bnp$ and $\Bns$. As a quantum group, or rather as a $C^*$-algebra given by some relations on the generators $u_{ij}$, our definition of $\Bns$ is just a new name for $\Bnp$ of \cite{BS09}. But as a category of noncrossing partitions our definition of $\Bnp$ coincides with $\Bnp$ of \cite{BS09}. This will be worked out in the sequel, see  Proposition \ref{PropKategorien}.

It might be confusing that we rename the former quantum group $\Bnp$ to $\Bns$ and that our new definition of a quantum group is named $\Bnp$ again, but the conceptual explanation for this is given by Propositions \ref{PropKonstrIsomI} and \ref{PropKonstrIsomII} and by Remark \ref{RemY}.
\end{rem}

The above seven quantum groups may be displayed in the following way (this diagram is dual to the one above):
\begin{align*}
 \Bn &\qquad\subset &\Bnp &\qquad\subset &\Bns &\qquad\subset &\On\\
\upsubset\; & &\upsubset\; & & & &\upsubset\;\\
 \Sn &\qquad\subset &\Snp & &\subset &  &\Hn
\end{align*}

We will now study their intertwiner spaces or rather their corresponding categories of noncrossing partitions. For this, we write $\langle p_1,\ldots,p_n\rangle$ for the (minimal) category generated by the partitions $p_1,\ldots, p_n\in NC$. This means, we consider $p_1,\ldots,p_n$ together with the unit partition $\idpart$ and the pair partition $\paarpart$ as a base case. Then we form the ``categorial hull'' by taking tensor products, compositions, involutions and rotations in the sense of Definition \ref{DefFullKateg}. 
The category $\langle\emptyset\rangle$ is the category of all noncrossing partitions generated by the unit partition and the pair partition.

Throughout this section, we will mainly be concerned with categories of noncrossing partitions, but the notation $\langle p_1,\ldots,p_n\rangle$ can also be used to describe the minimal full category generated by possibly crossing partitions $p_1,\ldots,p_n\in P$ in the sense of Definition \ref{DefFullKateg}. If all partitions $p_1,\ldots,p_n$ are noncrossing, then $\langle p_1,\ldots,p_n\rangle$ is automatically a category of noncrossing partitions.

In the sequel, some specified partitions play an important role.
Let $\singleton$ be the \emph{singleton partition} in $NC(0,1)$ consisting of a single point, not connected to any other point. We also denote it by $\shortmid$ (short bar) not to be confused with the unit partition $\idpart\in NC(1,1)$ (long bar). By $\downarrow$, we denote its adjoint in $NC(1,0)$.
Let $\vierpart$ be the \emph{four block partition} $\{1,2,3,4\}$ in $NC(0,4)$ connecting four points in a row. If $\vierpart$ is in a category $NC_X$ of noncrossing partitions, then -- by rotation -- the partition $\vierpartrot\;\in NC(2,2)$ is in $NC_X$, too, and vice versa. Finally, let $\legpart$ be the \emph{positioner partition} $\{1\}\{2,4\}\{3\}$ in $NC(0,4)$ connecting the second with the fourth point, and having a singleton on the first and on the third point, respectively. Note that the third point is \emph{not} connected to the pair and that $\shortmid$ denotes a singleton.  We try to keep the notation as simple as possible and therefore use the notation $\shortmid$ in $\legpart$ instead of $\singleton$. Let $\legpartrot\;\in NC(2,2)$ be a rotated version of $\legpart$. By $\singleton\otimes\singleton\;\in NC(0,2)$, we denote the tensor product of two singletons, by $\doublesingletonrot$ its rotated version in $NC(1,1)$.

The next lemma explains the correspondence between certain partitions $p$ (or rather their induced linear maps $T_p$ in the intertwiner spaces) and certain relations on the level of $C^*$-algebras.

\begin{lem}\label{LemNCRelCSternRel} 
Let $G$ be an orthogonal quantum group in the sense of Definition \ref{DefOrthQG}, let $A$ be its corresponding $C^*$-algebra, generated by self-adjoint elements $u_{ij}$ where $i,j=1,\ldots,n$, and denote by $C_G$ the corresponding intertwiner space.
 \begin{itemize}
  \item[(a)] For $p=\paarpart$, the map $T_p$ is in $C_G$ if and only if the $u_{ij}$ fulfill the relations $\sum_{k=1}^n u_{ik}u_{jk}=\delta_{ij}$ for all $i,j=1,\ldots,n$.

For $p=\;\baarpart$, the map $T_p$ is in $C_G$ if and only if the $u_{ij}$ fulfill the relations $\sum_{k=1}^n u_{ki}u_{kj}=\delta_{ij}$ for all $i,j=1,\ldots,n$.
  \item[(b)] For $p=\;\vierpartrot$, the map $T_p$ is in $C_G$ if and only if the $u_{ij}$ fulfill the relations $u_{ik}u_{jk}=u_{ki}u_{kj}=0$ whenever $i\neq j$,  for $i,j,k=1,\ldots,n$.
  \item[(c)] For $p=\;\doublesingletonrot$, the map $T_p$ is in $C_G$ if and only if the $u_{ij}$ fulfill the relations $\sum_{k=1}^n u_{ik}=\sum_{k=1}^n u_{kj}$ for all $i,j=1,\ldots,n$
  \item[(d)] For $p=\;\legpartrot$, the map $T_p$ is in $C_G$ if and only if the $u_{ij}$ fulfill the relations $u_{ij}(\sum_{k=1}^n u_{kl})=(\sum_{k=1}^n u_{mk})u_{ij}$ for all $i,j,l,m=1\ldots,n$.
  \item[(e)] For $p=\;\singleton$, the map $T_p$ is in $C_G$ if and only if the $u_{ij}$ fulfill the relations $\sum_{k=1}^n u_{ik}=1$ for all $i=1,\ldots,n$.

For $p=\;\downarrow$, the map $T_p$ is in $C_G$ if and only if the $u_{ij}$ fulfill the relations $\sum_{k=1}^n u_{kj}=1$ for all $j=1,\ldots,n$.
 \end{itemize}
\end{lem}
\Proof For $p\in NC(k,l)$, the linear map $T_p:(\C^n)^{\otimes k}\to (\C^n)^{\otimes l}$ gives rise to an element $T_p\otimes 1\in M_{n^k\times n^l}(A)$, for a fixed basis $e_1,\ldots e_n$ of $\C^n$. By $u^{\otimes k}$, we denote the $n^k\times n^k$ matrix $(u_{i_1j_1}\ldots u_{i_kj_k})$. Let $\xi=e_{i_1}\otimes \ldots\otimes e_{i_k}\in (\C^n)^{\otimes k}$. We apply $T_p\otimes 1$ respectively $u^{\otimes k}$ to $\xi\otimes 1$ by $(T_p\otimes 1)(\xi\otimes 1)=T_p\xi\otimes 1$ respectively  $u^{\otimes k}(\xi\otimes 1)=\sum_{\alpha_1,\ldots,\alpha_k=1}^n e_{\alpha_1}\otimes \ldots\otimes e_{\alpha_k}\otimes u_{\alpha_1i_1}\ldots u_{\alpha_ki_k}$.

(a) For $p=\paarpart\in NC(0,2)$, the linear map $T_p:\C\to(\C^n)^{\otimes 2}$ is given by \linebreak $T_p(1)=\sum_{i=1}^n e_i\otimes e_i$. Now, $T_p\in C_G$ if and only if $T_p(1)\otimes 1=u^{\otimes 2}(T_p(1)\otimes 1)$. The computation of $u^{\otimes 2}(T_p(1)\otimes 1)$ yields:
\[u^{\otimes 2}(T_p(1)\otimes 1)
=\sum_{k=1}^n u^{\otimes 2}(e_k\otimes e_k\otimes 1)
=\sum_{i,j=1}^ne_i\otimes e_j\otimes \left(\sum_{k=1}^n u_{i k}u_{jk}\right)\]
Comparison of the coefficients yields $\sum_{k=1}^n u_{ik}u_{jk}=\delta_{ij}$.
For $p=\;\baarpart\;\in NC(2,0)$, we have a similar result, using $T_p(e_i\otimes e_j)=\delta_{ij}$.

(b) If $p=\;\vierpartrot\;\in NC(2,2)$, the map $T_p$ is given by $T_p(e_i\otimes e_j)=\delta_{ij}e_i\otimes e_i$. We compute:
\begin{align*}
u^{\otimes 2}(T_p\otimes 1)(e_i\otimes e_j\otimes 1)
&=\sum_{k,l=1}^n\delta_{ij}e_k\otimes e_l\otimes u_{ki }u_{li}\\
(T_p\otimes 1)u^{\otimes 2}(e_i\otimes e_j\otimes 1)
&=\sum_{k=1}^ne_k\otimes e_k\otimes u_{ki }u_{kj}
\end{align*}
Comparison of the coefficients yields $u_{ki}u_{kj}=0$ if $i\neq j$. If $i=j$, then $u_{ki}u_{li}=0$ whenever $k\neq l$.

(c) For $p=\;\doublesingletonrot$, we have $T_p(e_i)=\sum_{j=1}^ne_j$. A direct computation yields:
\begin{align*}
u^{\otimes 1} (T_p\otimes 1)(e_i\otimes 1)
&=\sum_{k=1}^n\left(e_k\otimes\left(\sum_{j=1}^n u_{kj}\right)\right)\\
(T_p\otimes 1)u^{\otimes 1}(e_i\otimes 1)
&=\sum_{k=1}^n\left(e_k\otimes\left(\sum_{j=1}^n u_{ji}\right)\right)
\end{align*}
Thus, the sum $\sum_{j=1}^n u_{kj}$ is independent of $k$, whereas $\sum_{j=1}^n u_{ji}$ is independent of $i$, and they coincide if and only if $T_p\otimes 1$ and $u^{\otimes 1}$ commute.

(d) If $p=\;\legpartrot$, then $T_p(e_j\otimes e_l)=\sum_{k=1}^ne_k\otimes e_j$. Hence:
\begin{align*}
u^{\otimes 2} (T_p\otimes 1)(e_j\otimes e_l\otimes 1)
&=\sum_{m,i=1}^n\left(e_m\otimes e_i\otimes\left(\sum_{k=1}^n u_{mk}\right)u_{ij}\right)\\
(T_p\otimes 1)u^{\otimes 2} (e_j\otimes e_l\otimes 1)
&=\sum_{m,i=1}^n\left(e_m\otimes e_i\otimes u_{ij}\left(\sum_{k=1}^n u_{kl}\right)\right)
\end{align*}
 
(e) Let $p=\;\singleton\;\in NC(0,1)$. Then $T_p:\C\to \C^n$ is given by $T_p(1)=\sum_{i=1}^ne_i$ and we have $u^{\otimes 1}(T_p\otimes 1)(1\otimes 1)=\sum_{i=1}^ne_i\otimes\left(\sum_{k=1}^nu_{ik}\right)$.

If $p=\;\downarrow\;\in NC(1,0)$, then $T_p(e_j)=1$.
\qqed

We need a technical lemma about the positioner partition $\legpart$. 

In a partition $p\in P$, \emph{blocks} are the smallest subpartitions, i.e. any maximal subset of connected points in $p$ gives rise to a block.

\begin{lem}\label{LemLegPart} The following statements about categories of partitions and the positioner partition $\legpart$ hold true.
 \begin{itemize}
  \item[(a)] The positioner partition $\legpart$ is in the category $\langle\singleton\rangle$ generated by the singleton.
  \item[(b)] The tensor product $\singleton\otimes\singleton$ of two singletons is in the category $\langle\legpart\rangle$.
  \item[(c)] The positioner partition $\legpart$ is in $\langle\singleton\otimes\singleton,\vierpart\rangle$.
  \item[(d)] Let $NC_X$ be a category of noncrossing partitions (or a full category of partitions) containing the positioner partition $\legpart$, and let $p\in NC_X$ be a partition containing a singleton $\singleton$ as a block. Then, the partitions obtained from $p$ by shifting the singleton to any other position are all in $NC_X$.
 \end{itemize}
\end{lem}

\quad

\Proof (a) The partition $\legpart$ may be constructed by two singletons, a pair partition and the category operations (we compose $\paarpart\in NC(0,2)$ with the partition $\singleton\otimes\idpart\otimes\singleton\otimes\idpart\in NC(2,4)$):

\setlength{\unitlength}{0.5cm}
\begin{center}
\begin{picture}(9,4)
 \put(0,2.5){\uppartii{1}{2}{4}}
 \put(0,0){\upparti{1}{1}}
 \put(0,0){\upparti{2}{2}}
 \put(0,0){\upparti{1}{3}}
 \put(0,0){\upparti{2}{4}}
 \put(5.5,1.5){$=$}
 \put(6,1){\upparti{1}{1}}
 \put(6,1){\uppartii{2}{2}{4}}
 \put(6,1){\upparti{1}{3}}
 \end{picture}
\end{center}

(b) Compose the partition $\legpart$ with a pair partition (and two unit partitions) to obtain $\singleton\otimes\singleton$.

(c) Compose the partition $\singleton\otimes\singleton\otimes\vierpart$ with a pair to obtain the partition $\singleton\otimes\dreipart$ in $\langle\singleton\otimes\singleton,\vierpart\rangle$. 
We compose it with $\idpart^{\otimes 3}\otimes\singleton\otimes\dreipart\otimes\idpart$ and then with $\idpart\otimes\idpart\otimes\paarpart^*\otimes\idpart\otimes\paarpart^*\otimes\idpart$.

\setlength{\unitlength}{0.5cm}
\begin{center}
\begin{picture}(14,5)
 \put(0,2.5){\upparti{1}{1}}
 \put(0,2.5){\uppartiii{2}{2}{3}{8}}
 \put(0,2.5){\upparti{1}{4}}
 \put(0,2.5){\uppartiii{1}{5}{6}{7}}
 \put(0,0){\upparti{2}{1}}
 \put(0,0){\upparti{2}{2}}
 \put(0,3){\partii{1}{3}{4}}
 \put(0,0){\upparti{2}{5}}
 \put(0,3){\partii{1}{6}{7}}
 \put(0,0){\upparti{2}{8}}
 \put(9.5,2){$=$}
 \put(10,1.5){\upparti{1}{1}}
 \put(10,1.5){\uppartii{2}{2}{4}}
 \put(10,1.5){\upparti{1}{3}}
 \end{picture}
\end{center}

Hence, we infer $\legpart\in\langle\singleton\otimes\singleton,\vierpart\rangle$.

(d) Since $NC_X$ is a category, we also have $\legpartrot\in NC_X$ by rotation. 
We compose $\idpart\otimes\legpartrot$ with $\legpartrot\otimes\idpart$ and we deduce that $\singleton\otimes\idpart^{\otimes 2}\otimes\downarrow$ is in $NC_X$, too.

\setlength{\unitlength}{0.5cm}
\begin{center}
\begin{picture}(15,5)
 \put(0,2.5){\upparti{2}{2}}
 \put(0,2.5){\upparti{1}{3}}
 \put(0,2.5){\upparti{2}{4}}
 \put(0,5.5){\parti{1}{5}}
 \put(0,0){\upparti{1}{1}}
 \put(0,0){\upparti{2}{2}}
 \put(0,3){\parti{1}{3}}
 \put(0,0){\upparti{2}{4}}
 \put(6.5,2){$=$}
 \put(7,1.5){\upparti{1}{1}}
 \put(7,1.5){\upparti{2}{2}}
 \put(7,1.5){\upparti{2}{3}}
 \put(7,4.5){\parti{1}{4}}
 \put(12,2){$= \quad \singleton\otimes\;\idpart^{\otimes 2}\;\otimes\downarrow$}
 \end{picture}
\end{center}

By induction, we conclude that $\singleton\otimes\;\idpart^{\otimes k}\otimes\downarrow\;\in NC_X$ for all $k\in \N$, as well as their adjoints $\downarrow\otimes\;\idpart^{\otimes k}\otimes\singleton\;\in NC_X$. Using them, we can place any singleton in $p\in NC_X$ to any position. We did not make use of the fact that all partitions in $NC_X$ are noncrossing, thus the statement holds for arbitrary full categories of partitions, too.
\qqed

We will now describe the categories corresponding to the seven quantum groups of Lemma \ref{LemSiebenQG}. For the case of
 $\Bns$, we label the points of the partitions alternating by $\oplus\ominus\oplus\ominus\ldots\oplus\ominus$. Since there will be only partitions with an even number of blocks, we will have as many $\oplus$ as $\ominus$. For a general partition in $p\in NC(k,l)$, we start our labeling at the very right upper point with $\oplus$, and continue from right to  left in the upper row of points. We then proceed from left to right in the lower row of points, ending with $\ominus$. Thus, $p$ is labeled counterclockwise.

The next proposition refines the statement of \cite[th. 3.13]{BS09}.

\begin{prop}\label{PropKategorien}
 The following categories of noncrossing partitions correspond to the seven quantum groups of Lemma \ref{LemSiebenQG}.
\begin{itemize}
 \item[(1)] $\langle \emptyset\rangle=NC_2\subset NC$ is the category of all noncrossing pair partitions. It corresponds to $\On$.
 \item[(2)] $\langle \vierpart\rangle\subset NC$ is the category of all noncrossing partitions with blocks of even size (i.e. blocks with an even number of legs). It corresponds to $\Hn$.
 \item[(3)] $\langle \singleton\otimes\singleton,\vierpart\rangle\subset NC$ is the category of all noncrossing partitions with an even number of blocks of odd size (and an arbitrary number of blocks of even size). It corresponds to $\Snp$.
 \item[(4)] $\langle \singleton,\vierpart\rangle= NC$ is the category of all noncrossing partitions. It corresponds to $\Sn$.
 \item[(5)] $\langle \singleton\otimes\singleton\rangle\subset NC$ is the category of all noncrossing partitions with an arbitrary number of blocks of size two (pairs), each connecting one $\oplus$ with one $\ominus$, and an even number of blocks of size one (singletons). It corresponds to $\Bns$.
 \item[(6)] $\langle \legpart\rangle\subset NC$ is the category of all noncrossing partitions with an arbitrary number of blocks of size two and an even number of blocks of size one. It corresponds to $\Bnp$.
 \item[(7)] $\langle \singleton\rangle\subset NC$ is the category of all noncrossing partitions with blocks of size one or two. It corresponds to $\Bn$.
\end{itemize}
\end{prop}
\Proof (1) Clearly $NC_2$ is a category, since tensor products, compositions and involutions of noncrossing pairs are noncrossing pairs again. Furthermore, every noncrossing pair partition may be constructed out of the unit partition $\idpart$ and the pair partition $\paarpart$ using the category operations, hence $\langle \emptyset\rangle= NC_2$. By \cite[th. 3.13]{BS09} or \cite{BC07}, we know that $C_{\On}=\textnormal{span}\{T_p\;|\;p\in NC_2\}$. See also Lemma \ref{LemNCRelCSternRel}(a).

(2) The set of noncrossing partitions with blocks of even size is a category, since it is stable by the category operations (for the composition, note that the size of composed blocks is the addition of all its block sizes minus two times the relevant removed points). We denote it by $\mathcal C_2$. It contains the four block partition, thus also the category $\langle\vierpart\rangle$ generated by it. On the other hand, $\langle\vierpart\rangle$ contains every block partition (i.e. a partition consisting of a single block) of even size, by the following inductive argument ($n>1$):
\setlength{\unitlength}{0.5cm}
\newsavebox{\boxzwein}
\newsavebox{\boxzweinpluseins}
\newsavebox{\boxpaarundunit}
\savebox{\boxzwein}
{\begin{picture}(10,2)
\put(-0.1,0.5){\uppartiv{1}{1}{2}{5}{6}}
\put(-0.1,0.5){\uppartiv{1}{7}{8}{9}{10}}
\put(3.3,1.8){$2n$}
\put(3.3,0.5){$\ldots$}
\end{picture}}
\savebox{\boxzweinpluseins}
{\begin{picture}(6,2)
\put(-0.1,0.5){\uppartiv{1}{1}{2}{5}{6}}
\put(2.3,1.8){$2(n+1)$}
\put(3.3,0.5){$\ldots$}
\end{picture}}
\savebox{\boxpaarundunit}
{\begin{picture}(10,2)
\put(-0.1,0.5){\parti{2}{1}}
\put(-0.1,0.5){\parti{2}{2}}
\put(-0.1,0.5){\parti{2}{5}}
\put(-0.1,0.5){\partii{1}{6}{7}}
\put(-0.1,0.5){\parti{2}{8}}
\put(-0.1,0.5){\parti{2}{9}}
\put(-0.1,0.5){\parti{2}{10}}
\end{picture}}

\begin{center}
\begin{picture}(18,4)
 \put(0,2){\usebox{\boxzwein}}
 \put(0,2.5){\usebox{\boxpaarundunit}}
 \put(11.5,1.8){$=$}
 \put(12,1){\usebox{\boxzweinpluseins}}
 \end{picture}
\end{center}

Thus, using the category operations, every noncrossing partition with blocks of even size may be constructed in $\langle\vierpart\rangle$, hence $\mathcal C_2=\langle\vierpart\rangle$. By Lemma \ref{LemNCRelCSternRel}(b), the category $\langle\vierpart\rangle$ corresponds to $\Hn$ (note that $\vierpart\;\in NC_X$ if and only if $\vierpartrot\;\in NC_X$). (cf. also \cite[th. 3.13]{BS09} or \cite{BBC07})

(3) The set $\mathcal C_3$ of all noncrossing partitions with an even number of blocks of odd size is stable by the category operations. Indeed, if two blocks of odd size are getting connected by composition, the result is a block of even size. 
Thus, if $p$ and $q$ are two partitions in $\mathcal C_3$ with $m_p$ resp. $m_q$ blocks of odd size, then the number $m_{pq}$ of blocks of odd size in the composition of $p$ and $q$ is the sum of $m_p$ and $m_q$ minus an even number.
Therefore, $\mathcal C_3$ is a category containing $\langle \singleton\otimes\singleton,\vierpart\rangle$. 

Conversely, every partition $p\in \mathcal C_3$ is in $\langle \singleton\otimes\singleton,\vierpart\rangle$. To see this, replace each block $b_i$ of odd size $k_i$ in $p$ by a block $\tilde b_i$ of size $k_i+1$, respectively. The blocks of even size remain untouched. The resulting partition $p'$ is in $\langle \singleton\otimes\singleton,\vierpart\rangle$, it even is in $\langle \vierpart\rangle$ by (2). Hence, the partition $(\singleton\otimes\singleton)^{\otimes m}\otimes p'$ is in $\langle \singleton\otimes\singleton,\vierpart\rangle$, where $2m$ is the number of odd blocks of $p$. By Lemma \ref{LemLegPart}(c) and (d), we can shift a singleton besides each modified block $\tilde b_i$, respectively. Now, we use the pair partition to combine a singleton and a block $\tilde b_i$ to obtain $b_i$ back. Thus $p$ is in $\langle \singleton\otimes\singleton,\vierpart\rangle$.

By Lemma \ref{LemNCRelCSternRel}(b) and (c), $\mathcal C_3$ corresponds to $\Snp$, see also \cite[th. 3.13]{BS09}.

(4) By (2), $\langle \singleton,\vierpart\rangle$ contains every block partition of even size. Composing them with singletons $\singleton$ yield all block partitions of odd size. Now, by the category operations, all noncrossing partitions may be constructed, hence $\langle \singleton,\vierpart\rangle=NC$.
By Lemma \ref{LemNCRelCSternRel}(b) and (e), it corresponds to $\Sn$. (cf. also \cite[th. 3.13 or th. 3.8]{BS09})

(5) The set $\mathcal C_5$ of all noncrossing partitions with an arbitrary number of blocks of size two connecting $\oplus$ and $\ominus$, and an even number of blocks of size one is a category. To check the stability by tensor product, we swap the labels $\oplus$ and $\ominus$ in one of the partitions, if necessary. To see the stability under composition, recall that the counterclockwise labelling of a partition starts at the right upper point with $\oplus$ and ends at the right lower point with $\ominus$. Thus, the labelling of the middle row cancels out and we obtain a partition, which is in $\mathcal C_5$ again.
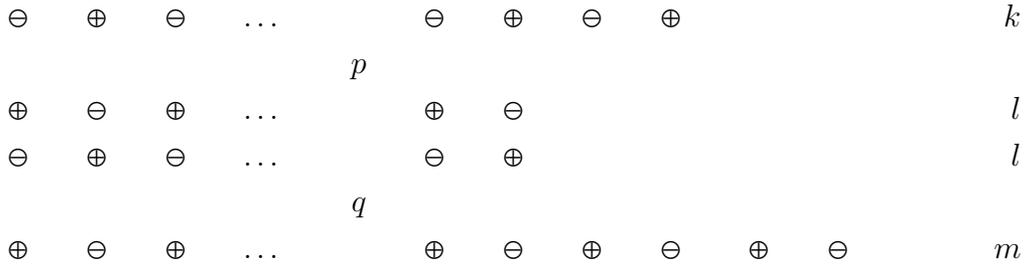
\begin{align*}
 \ominus &&\oplus &&\ominus &&\ldots &&\quad &&\ominus &&\oplus &&\ominus &&\oplus &&\quad && && &&k\\
&& && && && p\\
\oplus &&\ominus &&\oplus &&\ldots &&\quad &&\oplus &&\ominus && && && && && &&l\\
\ominus &&\oplus &&\ominus &&\ldots &&\quad  &&\ominus &&\oplus&& && && && && &&l\\
&& && && && q\\
 \oplus &&\ominus &&\oplus &&\ldots &&\quad  &&\oplus &&\ominus &&\oplus &&\ominus &&\oplus &&\ominus &&\quad &&m
\end{align*}
Since now $\mathcal C_5$ is a category containing $\singleton\otimes\singleton$, we have $\langle \singleton\otimes\singleton\rangle\subset \mathcal C_5$. To prove $\langle \singleton\otimes\singleton\rangle\supset \mathcal C_5$, we connect the singletons in a partition $p\in\mathcal C_5$ virtually to pairs in a noncrossing way. (For this, we may assume $p$ to be in $NC(0,l)$ by rotation, see Definition \ref{DefFullKateg}.) An example:
\setlength{\unitlength}{0.5cm}
\newsavebox{\dashlinehoch}
\newsavebox{\dashlinequer}
\newsavebox{\dashparti}
\newsavebox{\dashpartii}
\savebox{\dashlinehoch}
{\begin{picture}(1,1)
\put(0,0.25){\line(0,1){0.25}}
\put(0,0.75){\line(0,1){0.25}}
\end{picture}}
\savebox{\dashlinequer}
{\begin{picture}(1,1)
\put(0,0){\line(1,0){0.125}}
\put(0.375,0){\line(1,0){0.25}}
\put(0.875,0){\line(1,0){0.125}}
\end{picture}}
\savebox{\dashparti}
{\begin{picture}(3,1)
\put(0,0){\usebox{\dashlinehoch}}
\put(0,1){\usebox{\dashlinequer}}
\put(1,1){\usebox{\dashlinequer}}
\put(2,1){\usebox{\dashlinequer}}
\put(3,0){\usebox{\dashlinehoch}}
\end{picture}}
\savebox{\dashpartii}
{\begin{picture}(7,2)
\put(0,0){\usebox{\dashlinehoch}}
\put(0,1){\usebox{\dashlinehoch}}
\put(0,2){\usebox{\dashlinequer}}
\put(1,2){\usebox{\dashlinequer}}
\put(2,2){\usebox{\dashlinequer}}
\put(3,2){\usebox{\dashlinequer}}
\put(4,2){\usebox{\dashlinequer}}
\put(5,2){\usebox{\dashlinequer}}
\put(6,2){\usebox{\dashlinequer}}
\put(7,0){\usebox{\dashlinehoch}}
\put(7,1){\usebox{\dashlinehoch}}
\end{picture}}

\setlength{\unitlength}{0.5cm}
\begin{center}
\begin{picture}(14,4)
\put(-0.1,0){\upparti{1}{1}}
\put(-0.1,0){\uppartii{1}{2}{3}}
\put(-0.1,0){\upparti{1}{4}}
\put(-0.1,0){\upparti{1}{5}}
\put(-0.1,0){\uppartii{1}{6}{7}}
\put(-0.1,0){\uppartii{2}{8}{11}}
\put(-0.1,0){\uppartii{1}{9}{10}}
\put(-0.1,0){\upparti{1}{12}}
\put(-0.1,0){\uppartii{1}{13}{14}}
\put(1.15,1){\usebox{\dashparti}}
\put(5.15,1){\usebox{\dashpartii}}
\end{picture}
\end{center}
This cannot be done in a unique way, of course, since we might have a choice which singletons to connect. For instance, the four singletons in the preceding example may also be connected in the following way:
\setlength{\unitlength}{0.5cm}
\newsavebox{\dashpartiii}
\newsavebox{\dashpartiv}
\savebox{\dashpartiii}
{\begin{picture}(11,2)
\put(0,0){\usebox{\dashlinehoch}}
\put(0,1){\usebox{\dashlinehoch}}
\put(0,2){\usebox{\dashlinequer}}
\put(1,2){\usebox{\dashlinequer}}
\put(2,2){\usebox{\dashlinequer}}
\put(3,2){\usebox{\dashlinequer}}
\put(4,2){\usebox{\dashlinequer}}
\put(5,2){\usebox{\dashlinequer}}
\put(6,2){\usebox{\dashlinequer}}
\put(7,2){\usebox{\dashlinequer}}
\put(8,2){\usebox{\dashlinequer}}
\put(9,2){\usebox{\dashlinequer}}
\put(10,2){\usebox{\dashlinequer}}
\put(11,0){\usebox{\dashlinehoch}}
\put(11,1){\usebox{\dashlinehoch}}
\end{picture}}
\savebox{\dashpartiv}
{\begin{picture}(1,1)
\put(0,0){\usebox{\dashlinehoch}}
\put(0,1){\usebox{\dashlinequer}}
\put(1,0){\usebox{\dashlinehoch}}
\end{picture}}

\setlength{\unitlength}{0.5cm}
\begin{center}
\begin{picture}(14,4)
\put(-0.1,0){\upparti{1}{1}}
\put(-0.1,0){\uppartii{1}{2}{3}}
\put(-0.1,0){\upparti{1}{4}}
\put(-0.1,0){\upparti{1}{5}}
\put(-0.1,0){\uppartii{1}{6}{7}}
\put(-0.1,0){\uppartii{2}{8}{11}}
\put(-0.1,0){\uppartii{1}{9}{10}}
\put(-0.1,0){\upparti{1}{12}}
\put(-0.1,0){\uppartii{1}{13}{14}}
\put(1.15,1){\usebox{\dashpartiii}}
\put(4.15,1){\usebox{\dashpartiv}}
\end{picture}
\end{center}
But there is at least one possible virtual connection. We prove this by induction on the number $m$ of singletons in $p\in \mathcal C_5$. Let $m=2$ and assume the two singletons may not be connected in a noncrossing way. 
Thus, there is a pair which causes a crossing, if we connect the two singletons. Therefore, one singleton must be between the legs of this pair, whereas the other one is outside. To be more precise: If $a,b\in \{1,\ldots,l\}$, $a<b$ are the positions of the singletons and $x<y$ is the described pair, then $a<x<b<y$ or $x<a<y<b$.
\setlength{\unitlength}{0.5cm}
\newsavebox{\dashpartv}
\savebox{\dashpartv}
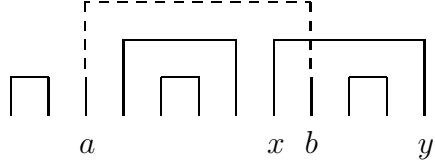
{\begin{picture}(6,2)
\put(0,0){\usebox{\dashlinehoch}}
\put(0,1){\usebox{\dashlinehoch}}
\put(0,2){\usebox{\dashlinequer}}
\put(1,2){\usebox{\dashlinequer}}
\put(2,2){\usebox{\dashlinequer}}
\put(3,2){\usebox{\dashlinequer}}
\put(4,2){\usebox{\dashlinequer}}
\put(5,2){\usebox{\dashlinequer}}
\put(6,0){\usebox{\dashlinehoch}}
\put(6,1){\usebox{\dashlinehoch}}
\end{picture}}

\setlength{\unitlength}{0.5cm}
\begin{center}
\begin{picture}(12,5)
\put(-0.1,1){\uppartii{1}{1}{2}}
\put(-0.1,1){\upparti{1}{3}}
\put(-0.1,1){\uppartii{2}{4}{7}}
\put(-0.1,1){\uppartii{1}{5}{6}}
\put(-0.1,1){\uppartii{2}{8}{12}}
\put(-0.1,1){\upparti{1}{9}}
\put(-0.1,1){\uppartii{1}{10}{11}}
\put(3.15,2){\usebox{\dashpartv}}
\put(3,0){$a$}
\put(8,0){$x$}
\put(9,0){$b$}
\put(12,0){$y$}
\end{picture}
\end{center}
We infer that there must be an odd number of points between $x$ and $y$ (note that the partition $p$ is noncrossing, if we do \emph{not} connect the singletons), thus $p$ does not follow the rule on $\oplus$ and $\ominus$, which is a contradiction.

For $m\geq4$, we denote the positions of the singletons by $i_1<i_2<\ldots<i_m$ and shift one singleton to position one, by rotation, thus $i_1=1$. 
If we may connect the singletons on $i_1$ and $i_2$ such that the resulting partition $p'$ is noncrossing, then $p'\in \mathcal C_5$, since the number of points between $i_1$ and $i_2$ is even. (Therefore, the pairing $\{i_1,i_2\}$ follows the rule on $\oplus$ and $\ominus$) Hence, we may apply the induction hypothesis. Otherwise, there is a pair connecting $x$ and $y$ such that $i_1<x<i_2<y$. Since $p$ follows the rule on $\oplus$ and $\ominus$, there must be an even number of singletons between $x$ and $y$, hence we have $i_1<x<i_2<i_{2k-1}<y<i_{2k}$ for a $k\geq 2$. We apply the induction hypothesis to the subpartition of $p$, connecting the points $\{x+1,\ldots,y-1\}$. We apply it a second time to its complement, i.e. the subpartition of $p$, connecting he remaining points of $\{1,\ldots,l\}$.

This finishes the proof that we may see the partitions in $\mathcal C_5$ as noncrossing pair partitions with two kinds of pairs: the usual ones ($\paarpart$) and the virtual ones, arising from two singletons. Thus, any partition in $\mathcal C_5$ may be constructed out of the pair partition $\paarpart$ and the ``virtual pair partition'' $\singleton\otimes\singleton$ by the category operations, hence $\langle \singleton\otimes\singleton\rangle\supset \mathcal C_5$. 

The correspondence of $\mathcal C_5$ and $\Bns$ is due to Lemma \ref{LemNCRelCSternRel}(c).

(6) The set $\mathcal C_6$ of all noncrossing partitions with an arbitrary number of blocks of size two and an even number of blocks of size one is a category, containing $\legpart$. On the other hand $\langle \legpart\rangle\supset \mathcal C_6$. To see this, let $p\in \mathcal C_6$ and denote by $p'$ its corresponding partition, if all singletons in $p$ are removed. Then $p'$ is in $NC_2\subset\langle\legpart\rangle$. By Lemma \ref{LemLegPart}(b), the partition $(\singleton\otimes\singleton)^{\otimes m}\otimes p'$ is in $\langle \legpart\rangle$, where $2m$ is the number of singletons in $p$. By Lemma \ref{LemLegPart}(d), we may shift the singletons such that we obtain $p$, thus $p$ is in $\langle \legpart\rangle$.
This category corresponds to $\Bnp$ by Lemma \ref{LemNCRelCSternRel}(d).

(7) The set of all noncrossing partitions with blocks of size one or two is a category containing $\singleton$, whereas on the other hand, any such partition is in  $\langle \singleton\rangle$. By Lemma \ref{LemNCRelCSternRel}(e), it corresponds to $\Bn$. (cf. also \cite[th. 3.13]{BS09})
\qqed

On the level of categories of noncrossing partitions, the seven quantum groups may be displayed as follows.
\begin{align*}
 \langle\singleton\rangle &\qquad\supset &\langle\legpart\rangle &\qquad\supset &\langle\singleton\otimes\singleton\rangle &\qquad\supset &\langle\emptyset\rangle=NC_2\\
\downsubset\; & &\downsubset\; & & & &\downsubset\;\\
\quad\\
 \langle\singleton, \vierpart\rangle=NC &\qquad\supset &\langle\singleton\otimes\singleton, \vierpart\rangle & &\supset &  &\langle\vierpart\rangle
\end{align*}

All categories but (5) of the preceding proposition have been described in \cite{BS09} (although not in terms of generating partitions). 
See also Remark \ref{RemX}.

In Proposition \ref{PropKategorien} we have seen that we can characterize the above seven categories by the structure of their blocks. To prove that there are exactly seven categories of noncrossing partitions, we need to study more on the relations between partitions and their blocks.

It is clear that if a noncrossing partition $p\in NC$ consists of blocks $b_1,\ldots,b_s$, which are all in a category $NC_X$ of noncrossing partitions, then $p$ is also in $NC_X$. This follows, since $p$ may be constructed out of its blocks by the category operations (cf. also \cite[lem. 3.15]{BS09}). But the crucial point will be, if the converse direction is true as well: If $p$ is in $NC_X$ and if $b\subset p$ is a block -- do we have $b\in NC_X$? 

Like before, we use the pair partition $\paarpart$ to erase the points between the legs of $b$. If there is an even number of points between two legs, we can erase all of them, otherwise we  end up with a singleton. By this, we see that the \emph{pair reduced version $p'$ of $p$ (with respect to $b$)} is in $NC_X$. It is (by rotation) of the following form:
\setlength{\unitlength}{0.5cm}
\begin{center}
\begin{picture}(22,5)
\put(-0.1,2){\upparti{1}{1}}
\put(-0.1,2){\uppartviii{2}{2}{5}{7}{10}{14}{17}{19}{22}}
\put(-0.1,2){\upparti{1}{6}}
\put(-0.1,2){\upparti{1}{11}}
\put(-0.1,2){\upparti{1}{18}}
\put(3.25,2){$\ldots$}
\put(2.2,1.7){$\underbrace{\qquad\qquad}$}
\put(3.25,0){$\beta_1$}
\put(8.25,2){$\ldots$}
\put(7.2,1.7){$\underbrace{\qquad\qquad}$}
\put(8.25,0){$\beta_2$}
\put(12.25,2){$\ldots$}
\put(15.25,2){$\ldots$}
\put(14.2,1.7){$\underbrace{\qquad\qquad}$}
\put(15.25,0){$\beta_{n-1}$}
\put(20.25,2){$\ldots$}
\put(19.2,1.7){$\underbrace{\qquad\qquad}$}
\put(20.25,0){$\beta_n$}
\put(11,4.25){$b$}
\put(-1,2.5){$p'\;=$}
\end{picture}
\end{center}

Here, the numbers $\beta_i$ of legs of $b$ between every two singletons are strictly positive respectively (otherwise we could remove two singletons using the pair partition).
To be more precise, we use the rotation to make sure that the very right point of $p$ is a point of $b$. Then, we compose $p$ with an appropriate partition of the form $\paarpart^{*\otimes\alpha_1}\otimes\idpart^{\otimes\gamma_1}\otimes\ldots\otimes\paarpart^{*\otimes\alpha_k}\otimes\idpart^{\otimes\gamma_k}$ such that all points of the resulting partition either are disconnected or belong to $b$. By this, we obtain the partition $p'$.

But we can improve this reduction of $p$.
The next lemma corresponds to lemma 3.15 of \cite{BS09}. We state it in a more detailed version. It is \emph{not} true for arbitrary categories of partitions $P_X\not\subset NC$.

\begin{lem}\label{LemDreiFuenfzehn}
 Let $NC_X$ be a category of noncrossing partitions and let $b\subset p$ be a block of a partition $p\in NC_X$.
\begin{itemize}
 \item[(a)] If $\singleton\;\in NC_X$, then $b\in NC_X$.
\item[(b)] If $\singleton\;\notin NC_X$ and $\singleton\otimes\singleton\;\notin NC_X$, then $b\in NC_X$.
\item[(c)] If $\singleton\;\notin NC_X$ and $\singleton\otimes\singleton\;\in NC_X$, then $\begin{cases} b\in NC_X & \textnormal{if } b \textnormal{ is of even length.} \\ 
\singleton\otimes b\in NC_X & \textnormal{if } b \textnormal{ is of odd length.}\end{cases}$
\end{itemize}
If $\singleton\;\notin NC_X$, then all partitions $p\in NC_X$ are of even length.
\end{lem}
\Proof (a) We can use the singleton to erase the singletons in the pair reduced version of $p$.

(b) Assume, that the pair reduced version $p'$ of $p$ does not equal $b$. We use the pair partition to reduce $p'$  either to the singleton $\singleton$ or to $\singleton\otimes\singleton$, which is a contradiction.

(c) The pair reduced version $p'$ of $p$ is of even length (otherwise we could construct the singleton, using the pair partition). If the positioner partition $\legpart$ is in $NC_X$, we simply shift the singletons of $p'$ side by side, by Lemma \ref{LemLegPart}(d), we erase them using the pair partition $\paarpart$ and the proof is complete. In the more general case, we can still do so, but we have to work to see this.

Let $n$ denote the number of singletons in $p'$. If $n=0$ or $n=1$, the statement holds. For $n\geq2$, we may assume (by rotation) that the number $\beta_1$ of legs between the first and the second singleton is minimal, i.e. $\beta_1\leq\beta_i$ for all $i$. Since all $\beta_i$ are strictly positive, we may use the pair partition $(n-2)$-times to obtain the following partition $p''$ in $NC_X$ (we compose $p'$ with $\idpart^{\otimes(\beta_1+\beta_2+1)}\otimes\paarpart^*\otimes\idpart^{\otimes(\beta_3-1)}\otimes\paarpart^*\otimes\ldots\otimes\paarpart^*\otimes\idpart^{\otimes\beta_n}$):
\setlength{\unitlength}{0.5cm}
\begin{center}
\begin{picture}(10,5)
\put(-0.1,2){\upparti{1}{1}}
\put(-0.1,2){\uppartiv{2}{2}{5}{7}{10}}
\put(-0.1,2){\upparti{1}{6}}
\put(3.25,2){$\ldots$}
\put(2.2,1.7){$\underbrace{\qquad\qquad}$}
\put(3.25,0){$\beta_1$}
\put(8.25,2){$\ldots$}
\put(7.2,1.7){$\underbrace{\qquad\qquad}$}
\put(8.25,0){$\beta'_2$}
\put(-1,2.5){$p''\;=$}
\end{picture}
\end{center}

Note that $\beta'_2\geq\beta_2\geq\beta_1$. Since now $p''$ is of even length, we have \linebreak $\beta'_2=\beta_1+2m$, for a positive integer $m\in \N_0$. Using the pair partition $m$-times, we see (by rotation) that the following partition $p'''\in NC(\beta_1+1,\beta_1+1)$ is in $NC_X$:
\setlength{\unitlength}{0.5cm}
\begin{center}
\begin{picture}(6,8)
\put(-0.1,2){\upparti{1}{1}}
\put(-0.1,2){\uppartii{2}{2}{5}}
\put(3.75,4){\line(0,1){1}}
\put(-0.1,8){\partii{2}{2}{5}}
\put(-0.1,8){\parti{1}{6}}
\put(3.25,7){$\ldots$}
\put(3.25,2){$\ldots$}
\put(2.2,1.7){$\underbrace{\qquad\qquad}$}
\put(3.25,0){$\beta_1$}
\put(-2,4){$p'''\;=$}
\end{picture}
\end{center}

Finally, composing the partition $p'''$ with $p'$ effects, that we may shift the second singleton of $p'$ besides the first, and we erase them. Iterating this procedure yields the result.\qqed

\begin{thm}\label{SatzGenauSieben}
 There are exactly seven categories of noncrossing partitions, namely those of Proposition \ref{PropKategorien}.
\end{thm}
\Proof Let $NC_X$ be a category of noncrossing partitions.

\emph{Case 1.} Let $\singleton\;\in NC_X$ and $\vierpart\;\in NC_X$. By Proposition \ref{PropKategorien} we have that \linebreak $NC=\langle\singleton, \vierpart\rangle\subset NC_X\subset NC$.

\emph{Case 2.} Let $\singleton\;\in NC_X$ and $\vierpart\;\notin NC_X$, thus $\langle\singleton\rangle\subset NC_X$. By Lemma \ref{LemDreiFuenfzehn}, all blocks of partitions in $NC_X$ are in $NC_X$ and hence they are of length one or two (a block of bigger length could be shortened to a block of length three using the singleton, which would imply that the four block was in $NC_X$, see Section \ref{SectConcept}.)
Thus, $NC_X\subset\langle\singleton\rangle$ by Proposition \ref{PropKategorien}.

\emph{Case 3.} Let $\singleton\;\notin NC_X, \singleton\otimes\singleton\;\notin NC_X$ and $\vierpart\;\in NC_X$. By Lemma \ref{LemDreiFuenfzehn}, all blocks of partitions in $NC_X$ are of even length and in $NC_X$. By Proposition \ref{PropKategorien}, we conclude $NC_X=\langle \vierpart\rangle$.

\emph{Case 4.} Let $\singleton\;\notin NC_X, \singleton\otimes\singleton\;\notin NC_X$ and $\vierpart\;\notin NC_X$. Again, all blocks of partitions in $NC_X$ are of even size. But since $\vierpart\;\notin NC_X$, we infer that $NC_X\subset NC_2$. Now, $NC_2$ is contained in any category of noncrossing partitions, thus $NC_X= NC_2$.

\emph{Case 5.} Let $\singleton\;\notin NC_X, \singleton\otimes\singleton\;\in NC_X$ and $\vierpart\;\in NC_X$. Since $\singleton\;\notin NC_X$, all partitions in $NC_X$ are of even length. Thus, for any partition in $NC_X$ the number of blocks of odd size is even. We deduce $NC_X=\langle \singleton\otimes\singleton, \vierpart\rangle$, using Proposition \ref{PropKategorien}.

\emph{Case 6.} Let $\singleton\;\notin NC_X, \singleton\otimes\singleton\;\in NC_X, \vierpart\;\notin NC_X$ and $\legpart\;\in NC_X$. Let $p\in NC_X$ and let $b\subset p$ be a block of $p$. If $b$ is of even length, then $b\in NC_X$, by Lemma \ref{LemDreiFuenfzehn}. But since $\vierpart\;\notin NC_X$, $b$ must be of length two.
If $b$ is of odd length, then $\singleton\otimes b\in NC_X$, again by Lemma \ref{LemDreiFuenfzehn}. Thus, $b$ must be a singleton. Otherwise we could construct the four block, using the pair partition. We infer $NC_X\subset\langle\legpart\rangle$ by Proposition \ref{PropKategorien} and we conclude $NC_X=\langle\legpart\rangle$.

\emph{Case 7.} Let $\singleton\;\notin NC_X, \singleton\otimes\singleton\;\in NC_X, \vierpart\;\notin NC_X$ and $\legpart\;\notin NC_X$. Again, we have $NC_X\subset\langle\legpart\rangle$ by the above argumentation, but only $\langle\singleton\otimes\singleton\rangle\subset NC_X$. To see $\langle\singleton\otimes\singleton\rangle\supset NC_X$, let $p\in NC_X$ and assume that $p$ does not follow the rule on $\oplus$ and $\ominus$ as described in Proposition \ref{PropKategorien}. Thus, there is a pair with an odd number of points between its legs. By rotation, we may assume $p$ to be of the following form:
\setlength{\unitlength}{0.5cm}
\begin{center}
\begin{picture}(14,4)
\put(2,2){$p\;=$}
\put(1.9,2){\uppartii{1}{6}{11}}
\put(5.25,2){$\ldots$}
\put(4.2,1.7){$\underbrace{\qquad\qquad}$}
\put(0,0){odd number of points}
\put(10.25,2){$\ldots$}
\put(9.2,1.7){$\underbrace{\qquad\qquad}$}
\put(9,0){odd number of points}
\end{picture}
\end{center}

Using the pair partition we deduce, that $\legpart$ is in $NC_X$, which is a contradiction.
\qqed

\begin{cor}\label{GenauSieben}
 There are exactly seven free easy quantum groups (also called free orthogonal quantum groups), namely those of Lemma \ref{LemSiebenQG}.
\end{cor}

For the comparison of 
Corollary \ref{GenauSieben} with \cite[th. 3.16]{BS09} or rather with \cite[th. 1.4]{BCS11}, 
cf. Remark \ref{RemX}.

For the classical case, there are only six easy (quantum) groups. This is, we consider an easy quantum group $G$ with $S_n\subset G\subset O_n$. In fact, the associated $C^*$-algebra $A_{G}$ is generated by \emph{commuting} generators $u_{ij}$, thus $G$ is a classical group with $A_G=\mathcal C(G)$. The intertwiner space of an easy group is of the form \linebreak $C_{G}=\textnormal{span}\{T_p\;|\; p\in P_X\}$
where $P_X$ is a full category of partitions (cf. Definition \ref{DefEasyQG}).

Let $G$ be any orthogonal quantum group and let $p=\;\crosspart\;\in P(2,2)$ be the \emph{crossing partition} (also called the \emph{symmetry partition}) connecting two upper and two lower points by pairs in a crossing way. Then, the map $T_p$ is in the intertwiner space $C_G$ if and only if the generators $u_{ij}$ of $A_G$ commute (this may be proven analogously to Lemma \ref{LemNCRelCSternRel}). Thus, the classification of easy groups reduces to the classification of all full categories of partitions containing the crossing partition $\crosspart$.

Analogous to Definition \ref{DefSieben}, we define six orthogonal quantum groups $O_n$, $H_n$, $S_n'$, $S_n$, $B_n'$ and $B_n$ by their $C^*$-algebras with generators $u_{ij}$, $i,j=1,\ldots,n$ and exactly the same relations as in Definition \ref{DefSieben}. In addition, we require the $u_{ij}$ to commute mutually. Thus, the resulting quantum groups are in fact groups, listed by Banica and Speicher in \cite[prop. 2.4]{BS09}. Note that a system of commuting generators $u_{ij}$ which is bistochastic$^\#$ is automatically bistochastic'. Hence, the classical versions of $\Bnp$ and $\Bns$ coincide. (See also Remark \ref{RemY}.)

The full categories of partitions corresponding to these six easy groups are described in \cite[prop. 2.5]{BS09}. They are given by the categories (1)-(4), (6) and (7) as described in Proposition \ref{PropKategorien} by adding the crossing partition $\crosspart$ to the set of generators of each category. 
The verbal description of these categories is given by erasing the word ``noncrossing'' in (1)-(4), (6) and (7).
On the level of categories, we see that $\langle\singleton\otimes\singleton,\crosspart\;\rangle=\langle\legpart,\crosspart\;\rangle$, since $\legpart$ may be constructed by $\singleton\otimes\singleton\otimes\;\paarpart$ composed with the crossing partition $\crosspart$. This explains why the classical group versions of $\Bnp$ and $\Bns$ coincide.

By \cite[th. 2.8]{BS09}, there are exactly six easy groups, namely the above ones. The proof of Theorem \ref{SatzGenauSieben} (of our article) carries over to the  situation of categories containing the crossing partition $\crosspart$ (the cases 6 and 7 coincide), since Lemma \ref{LemDreiFuenfzehn} also holds true for full categories $P_X$ of partitions containing the crossing partition.
(If $b\subset p$ is a block of a partition $p\in P_X$, we have $b\otimes q\in P_X$ for some $q\in P$ using the crossing partition $\crosspart\;\in P_X$.)
We summarize the group case:

\begin{thm}[\cite{BS09}] \label{SatzGenauSechs}
 There are exactly six full categories of partitions containing the crossing partition $\crosspart$.
\end{thm}

\begin{cor}[\cite{BS09}]
 There are exactly six easy groups: $O_n$, $H_n$, $S_n'$, $S_n$, $B_n'$ and $B_n$.
\end{cor}

\section{The classification of the half-liberated easy quantum groups}\label{Sect4}

In order to classify all easy quantum groups in the sense of Definition \ref{DefEasyQG}, we have to study all possible full categories of partitions $P_X$ (recall Definition \ref{DefFullKateg}). In \cite[def. 6.6.]{BS09} and \cite[def. 2.2, th. 2.4]{BCS10}, further examples of easy quantum groups were defined, the so-called half-liberated quantum groups.

\begin{defn}
 An easy quantum group is called \emph{half-liberated}, if its category of partitions contains the \emph{half-liberating partition} $\halflibpart$ but not the crossing partition $\crosspart$. The half-liberating partition $\halflibpart\;\in P(3,3)$ is given by the partition $\{1,3'\}\{2,2'\}\{3,1'\}$ connecting three upper points $1,2,3$ (labeled from left to right) with three lower points $1',2',3'$ (labeled from left to right) such that $1$ and $3'$ are connected, $2$ and $2'$, and finally $3$ and $1'$.
\end{defn}
 
We will recall the definition and the categories of the two half-liberated easy quantum groups $\Onhl$ and $\Hnhl$ introduced in \cite{BS09} and \cite{BCS10}, and we will study a third one, $\Bnhl$. Also, we will reformulate the statement of \cite{BCS10} about the number of nonhyperoctahedral easy quantum groups.

The following half-liberated versions of free easy quantum groups are obtained by taking the quotient  by the relations $u_{ij}u_{kl}u_{pq}=u_{pq}u_{kl}u_{ij}$ for all $i,j,k,l,p,q$. In \cite[lem. 2.3]{BCS10} it is shown, that these relations correspond exactly to the half-liberating partition $\halflibpart$ in the sense of Lemma \ref{LemNCRelCSternRel}.

\begin{defn}\label{DefHalfLibQG}
 Let $n\in \N$ and let $G^+$ be a free easy quantum group. Denote by $\langle abc=cba\rangle$ the ideal in $A_{G^+}$ generated by the relations  $u_{ij}u_{kl}u_{pq}=u_{pq}u_{kl}u_{ij}$ for all $i,j,k,l,p,q$. We put:
\begin{itemize}
\item[(1)] $\AOnhl:= \AOn / \langle abc=cba\rangle$
\item[(2)] $\AHnhl:= \AHn / \langle abc=cba\rangle$
\item[(3)] $\ABnhl:= \ABns / \langle abc=cba\rangle$
\end{itemize}
\end{defn}

\begin{lem}\label{LemHalfLibQG}
 The three $C^*$-algebras of Definition \ref{DefHalfLibQG} give rise to orthogonal quantum groups $\Onhl$, $\Hnhl$ and $\Bnhl$.
\end{lem}
\Proof We only have to check that the $abc=cba$ relations pass through the maps $\Delta$, $\epsilon$ and $S$ of Definition \ref{DefOrthQG}.
 \qqed

\begin{rem}
 The half-liberated versions $\Onhl$ and $\Hnhl$ have been defined in \cite{BS09} and \cite{BCS10}. The quantum group $\Bnhl$ is new.

There are no half-liberated versions of $\Bn, \Bnp, \Sn$ and $\Snp$ since the quotients by  $\langle abc=cba\rangle$ are commutative $C^*$-algebras, thus $B_n^*=B_n$ etc. Indeed, denote by $r=\sum_{k=1}^nu_{ik}=\sum_{k=1}^nu_{kj}$ the symmetry in $\ABnp / \langle abc=cba\rangle$ (cf. Remark \ref{BemUeberR}). Now, $r$ commutes with all $u_{ij}$, thus:
\[u_{ij}u_{pq}=u_{ij}u_{pq}r^2=u_{ij}ru_{pq}r=\left(\sum_{l=1}^nu_{ij}u_{kl}u_{pq}\right) r=\left(\sum_{l=1}^nu_{pq}u_{kl}u_{ij}\right) r=u_{pq}u_{ij}\]
Since $A_{G^+} / \langle abc=cba\rangle$ is a quotient of $\ABnp / \langle abc=cba\rangle$, if $G^+=\Bn, \Sn$ or $\Snp$, there are no half-liberated versions of these easy quantum groups.
\end{rem}

Now, we describe the categories corresponding to these half-liberated quantum groups. The descriptions for $\Onhl$ and $\Hnhl$ are taken from \cite[th. 2.4]{BCS10}. By $P$ we denote the class of all noncrossing partitions. Again, we label the points of the partitions alternating by $\oplus \ominus \oplus \ominus \ldots \oplus \ominus$.
We observe that the half-liberating partition $\halflibpart$ acts like the crossing partition $\crosspart$ on every two points labeled by $\oplus$, and likewise on every two labeled by $\ominus$.

\begin{prop}\label{HalfLibKateg}
 The following full categories of partitions correspond to the three quantum groups of Lemma \ref{LemHalfLibQG} which hence are easy quantum groups.
\begin{itemize}
 \item[(1)] $\langle \halflibpart\rangle\subset P$ is the category of all pair partitions, each connecting one $\oplus$ with one $\ominus$. It corresponds to $\Onhl$.
 \item[(2)] $\langle \halflibpart, \vierpart\rangle\subset P$ is the category of all partitions with blocks of even size, each block consisting half of points labeled by $\oplus$ and half of those labeled by $\ominus$. It corresponds to $\Hnhl$.
 \item[(3)] $\langle \halflibpart, \singleton\otimes\singleton\rangle\subset P$ is the category of partitions with an arbitrary number of blocks of size two (pairs), each connecting one $\oplus$ with one $\ominus$, and an even number of blocks of size one (singletons). It corresponds to $\Bnhl$.
\end{itemize}
\end{prop}
\Proof (1) The set $\mathcal E_1$ of all pair partitions, each connecting one $\oplus$ with one $\ominus$, is a category (see also the proof of Proposition \ref{PropKategorien}(5)) containing $\halflibpart$, thus $\langle \halflibpart\rangle\subset \mathcal E_1$. Conversely, we obtain any partition of $\mathcal E_1$ by using the partition $\halflibpart$  to place the legs on the $\ominus$-points of $\paarpart^{\otimes k}$ in the right way. The correspondence to $\Onhl$ is due to \cite[lem. 2.3]{BCS10}.

(2) The set $\mathcal E_2$ of all partitions with blocks of even size, each block consisting half of points labeled by $\oplus$ and half of those labeled by $\ominus$, is a category containing $\halflibpart$ and $\vierpart$. On the other hand, we may construct any partition of $\mathcal E_2$ using $\halflibpart$ and $\vierpart$ and the category operations, analogous to the procedure described in the proof of (1). We use Lemma \ref{LemNCRelCSternRel}(b) and \cite[lem. 2.3]{BCS10} to deduce that $\mathcal E_2$ corresponds to $\Hnhl$.

(3) Again, the described set is a category whose partitions may be constructed applying $\halflibpart$ to tensor products of $\paarpart$ and $\singleton\otimes\singleton$. The correspondence to $\Bns$ is due to Lemma \ref{LemNCRelCSternRel}(c). \qqed

Next, we reformulate a result by Banica, Curran and Speicher (\cite{BCS10})
solving the classification of easy quantum groups in the nonhyperoctahedral case (which is defined below).
The following few lemmas are a preparation for this.

In the first lemma we prove that there are not only no half-liberated versions of $\Bn$, $\Bnp$, $\Sn$ and $\Snp$, but in fact there are no easy quantum groups at all in between $G$ and $G^+$, if $G=B_n$, $B_n'$, $S_n$ or $S_n'$.

\begin{lem}\label{LemmaA}
 Let $P_X\not\subset NC$ be a full category of partitions, containing the positioner partition $\legpart$. Then the crossing partition $\crosspart$ is in $P_X$.
\end{lem}
\Proof $P_X$ contains a partition $p_0\in P\backslash NC$ with a crossing.
By rotation it is of the following form, the block on $\alpha$ and $\gamma$ (and potentially on further points) causing a crossing with the block on $\beta$ and $\delta$ (and potentially on further points):
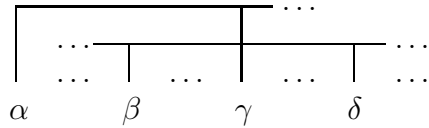
\setlength{\unitlength}{0.5cm}
\begin{center}
\begin{picture}(13,4)
\put(-0.1,1){\uppartii{2}{1}{7}}
\put(7,3){\line(1,0){1}}
\put(8.25,2.9){$\ldots$}
\put(-0.1,1){\uppartii{1}{4}{10}}
\put(3.25,2){\line(1,0){1}}
\put(2.25,1.9){$\ldots$}
\put(10,2){\line(1,0){1}}
\put(11.25,1.9){$\ldots$}
\put(2.25,1){$\ldots$}
\put(5.25,1){$\ldots$}
\put(8.25,1){$\ldots$}
\put(11.25,1){$\ldots$}
\put(1,0){$\alpha$}
\put(4,0){$\beta$}
\put(7,0){$\gamma$}
\put(10,0){$\delta$}
\end{picture}
\end{center}

We apply the rotated positioner partition $\legpartrot$ or rather $\singleton\otimes\idpart^{\otimes k}\otimes\downarrow$ (see also Lemma \ref{LemLegPart}(d), which is also true for full categories) to all points except $\alpha, \beta, \gamma$ and $\delta$ such that we obtain a partition consisting of two crossing pairs and some singletons. To be more precise, let $p_0\in P(0,l)$ be as above. We compose $p_0$ iterative with $\singleton\otimes\idpart^{\otimes k}\otimes\downarrow\otimes\idpart^{\otimes(l-k-1)}$ for all $k\in \{1,\ldots,l-1\}\backslash\{\beta-1,\gamma-1,\delta-1\}$ and we obtain:
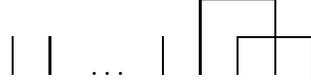
\setlength{\unitlength}{0.5cm}
\begin{center}
\begin{picture}(10,3)
\put(-0.1,0){\upparti{1}{1}}
\put(-0.1,0){\upparti{1}{2}}
\put(3.25,0){$\ldots$}
\put(-0.1,0){\upparti{1}{5}}
\put(-0.1,0){\uppartii{2}{6}{8}}
\put(-0.1,0){\uppartii{1}{7}{9}}
\end{picture}
\end{center}

We erase the singletons using the pair partition $\paarpart$, and by rotation we obtain the crossing partition $\crosspart$. (If the number of singletons in the above partition is odd, then the singleton $\singleton$ is in $P_X$ and we can use it to obtain the crossing partition from the above partition.) \qqed

\begin{cor}\label{CorDoubleSingl}
 Let $P_X$ be a full category of partitions with $P_X\not\subset NC$, $\vierpart\;\in P_X$ and $\crosspart\;\notin P_X$. Then $\singleton\otimes\singleton\;\notin P_X$.
\end{cor}
\Proof Assume $\singleton\otimes\singleton\;\in P_X$. By Lemma \ref{LemLegPart}(c), we then have $\legpart\;\in P_X$. By Lemma \ref{LemmaA}, we infer $\crosspart\;\in P_X$ which is a contradiction.\qqed

The next lemma is just a technical lemma used in Lemma \ref{LemmaC}.

\begin{lem}\label{LemmaB}
 Let $P_X$ be a full category of partitions containing the partition\linebreak $\singleton\otimes\crosspart\otimes\singleton$. Then the half-liberating partition $\halflibpart$ is in $P_X$.
\end{lem}
\Proof We use the pair partition $\paarpart$ to infer that $\singleton\otimes\singleton$ is in $P_X$. Thus, we can construct the partition $\halflibpart$ in $P_X$ in the following way using three times the partition $\singleton\otimes\crosspart\otimes\singleton$ in the rotated versions $\singleton\otimes\crosspart\otimes\downarrow$ and $\downarrow\otimes\crosspart\otimes\singleton$.
\setlength{\unitlength}{0.5cm}
\newsavebox{\kreuzer}
\savebox{\kreuzer}
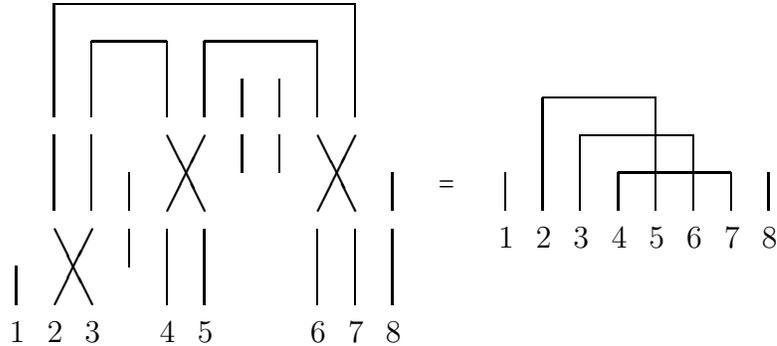
{\begin{picture}(1,2)
\thicklines
\put(0,2){\line(1,-2){1}}
\put(0,0){\line(1,2){1}}
\end{picture}}

\setlength{\unitlength}{0.5cm}
\begin{center}
\begin{picture}(19,10)
\put(0,6){\uppartii{3}{1}{9}}
\put(0,6){\uppartii{2}{2}{4}}
\put(0,6){\uppartii{2}{5}{8}}
\put(0,6){\upparti{1}{6}}
\put(0,6){\upparti{1}{7}}
\put(0,3.5){\upparti{2}{1}}
\put(0,3.5){\upparti{2}{2}}
\put(0,3.5){\upparti{1}{3}}
\put(4.3,3.5){\usebox{\kreuzer}}
\put(0,4.5){\upparti{1}{6}}
\put(0,4.5){\upparti{1}{7}}
\put(8.3,3.5){\usebox{\kreuzer}}
\put(0,3.5){\upparti{1}{10}}
\put(-1,1){\upparti{1}{1}}
\put(1.3,1){\usebox{\kreuzer}}
\put(0,2){\upparti{1}{3}}
\put(0,1){\upparti{2}{4}}
\put(0,1){\upparti{2}{5}}
\put(0,1){\upparti{2}{8}}
\put(0,1){\upparti{2}{9}}
\put(0,1){\upparti{2}{10}}
\put(0.1,0){1}
\put(1.1,0){2}
\put(2.1,0){3}
\put(4.1,0){4}
\put(5.1,0){5}
\put(8.1,0){6}
\put(9.1,0){7}
\put(10.1,0){8}
\put(11.5,4){$=$}
\put(12,3.5){\upparti{1}{1}}
\put(12,3.5){\uppartii{3}{2}{5}}
\put(12,3.5){\uppartii{2}{3}{6}}
\put(12,3.5){\uppartii{1}{4}{7}}
\put(12,3.5){\upparti{1}{8}}
\put(13.1,2.5){1}
\put(14.1,2.5){2}
\put(15.1,2.5){3}
\put(16.1,2.5){4}
\put(17.1,2.5){5}
\put(18.1,2.5){6}
\put(19.1,2.5){7}
\put(20.1,2.5){8}
\end{picture}
\end{center}
Using rotation and the pair partition $\paarpart$ we obtain $\halflibpart\;\in P_X$.\qqed

Next, we show that there are no easy quantum groups in between $O_n$ and $\On$ or $B_n'$ and $\Bns$ respectively apart from the half-liberated versions $\Onhl$ and $\Bnhl$.

\begin{lem}\label{LemmaC}
 Let $P_X$ be a full category of partitions, such that all partitions have blocks of length at most two. If $P_X\not\subset NC$, then either $\crosspart\;\in P_X$ or $\halflibpart\;\in P_X$.
\end{lem}
\Proof There is a partition $p_0\in P_X$ containing two crossing pairs as blocks.
By rotation and using the pair partition $\paarpart$, we may construct the following partition $p_0'$ from it. Here, $\alpha, \beta, \gamma$ and $\delta$ are either singletons, or connected by pairs, or zero, respectively:

\setlength{\unitlength}{0.5cm}
\begin{center}
\begin{picture}(8,4)
\put(-0.5,1){$p_0'\;\; =$}
\put(1.9,1.1){\uppartii{2}{1}{5}}
\put(1.9,1.1){\uppartii{1}{3}{7}}
\put(2,1){$\bullet$}
\put(4,1){$\bullet$}
\put(6,1){$\bullet$}
\put(8,1){$\bullet$}
\put(2,0){$\alpha$}
\put(4,0){$\beta$}
\put(6,0){$\gamma$}
\put(8,0){$\delta$}
\end{picture}
\end{center}

If $P_X$ contains a partition of odd length, then the singleton $\singleton$ is in $P_X$ (using the pair partition $\paarpart$ to erase all points but one), and we obtain $\crosspart\;\in P_X$ from $p_0'$. Otherwise, we have to distinguish three cases for the length $l$ of $p_0'$.

\emph{Case 1.} If $l=4$, then $p_0'$ is the crossing partition $\crosspart$ in a rotated form.

\emph{Case 2a.} If $l=6$ and two of the points $\alpha, \beta, \gamma$ and $\delta$ are singletons, we check all six cases for the positions of the singletons. Using the pair partition, we either obtain the positioner partition $\legpart$ or the partition $\singleton\otimes\crosspart\otimes\singleton$ in a rotated form, respectively. By Lemma \ref{LemmaA} and \ref{LemmaB}, we deduce $\crosspart\;\in P_X$ or $\halflibpart\;\in P_X$.

\emph{Case 2b.} If $l=6$ and $p_0'$ consists of three pairs, we again check all six cases for the position of the third pair and infer $\crosspart\;\in P_X$ or $\halflibpart\;\in P_X$.

\emph{Case 3.} If $l=8$ and $p_0'$ consists of two pairs and four singletons, we check that $\legpart\in P_X$ (thus $\crosspart\;\in P_X$ by Lemma \ref{LemmaA}). Likewise, if $p_0'$ consists of three pairs and two singletons (six cases). If $p_0'$ consists of four pairs, in all of these three cases we deduce $\crosspart\;\in P_X$. \qqed

\begin{lem}\label{LemmaD}
 Let $P_X$ be a full category of partitions containing a partition $p\in P_X$ such that there is a block $b\subset p$ of length at least three. Then the four block partition $\vierpart$ is in $P_X$.
\end{lem}
\Proof By rotation and using the pair partition $\paarpart$, we infer that the following partition $p'$ is in $P_X$. Here, $\alpha, \beta$ and $\gamma$ are either singletons, or connected by any partitions ($\gamma$ might even be connected to any of the points of the three block), or zero. Let $p''$ be the partition obtained from $(p')^*$ after rotation.

\setlength{\unitlength}{0.5cm}
\begin{center}
\begin{picture}(16,3)
\put(-0.5,1){$p'\;\; =$}
\put(0.9,1.1){\uppartiii{1}{1}{3}{5}}
\put(7.15,1.1){\usebox{\dashlinehoch}}
\put(6.15,2.1){\usebox{\dashlinequer}}
\put(3,1){$\bullet$}
\put(5,1){$\bullet$}
\put(7,1){$\bullet$}
\put(3,0){$\alpha$}
\put(5,0){$\beta$}
\put(7,0){$\gamma$}
\put(9.5,1){$p''\;\; =$}
\put(11.9,1.1){\uppartiii{1}{1}{3}{5}}
\put(12.15,1.1){\usebox{\dashlinehoch}}
\put(12.15,2.1){\usebox{\dashlinequer}}
\put(12,1){$\bullet$}
\put(14,1){$\bullet$}
\put(16,1){$\bullet$}
\put(12,0){$\gamma$}
\put(14,0){$\beta$}
\put(16,0){$\alpha$}
\end{picture}
\end{center}
Applying the pair partition $\paarpart$ to $p'\otimes p''$ yields that the following partition is in $P_X$.
(Compose $p'\otimes p''$ iterative with $\idpart^{\otimes a_i}\otimes\paarpart^*\otimes\idpart^{\otimes a_i}$ with appropriate $a_i\in \N$.)

\setlength{\unitlength}{0.5cm}
\begin{center}
\begin{picture}(6,3)
\put(0,1){\uppartiv{1}{1}{3}{4}{6}}
\put(2,1){$\bullet$}
\put(5,1){$\bullet$}
\put(2,0){$\alpha$}
\put(5,0){$\alpha$}
\end{picture}
\end{center}

We apply again the pair partition $\paarpart$ to the tensor product of this partition with its copy and we obtain the four block partition $\vierpart$ in $P_X$.
\qqed

We are now prepared to count all full categories of partitions that do \emph{not} contain the four block partition $\vierpart$.

\begin{prop}\label{PropNonhyperoct}
 Let $P_X$ be a full category of partitions.
\begin{itemize}
 \item[(a)] If $P_X\subset NC$, then $P_X$ is one of the seven cases of Theorem \ref{SatzGenauSieben} or rather of Proposition \ref{PropKategorien}.
 \item[(b)] If $P_X\not\subset NC$, then:
     \begin{itemize}
      \item[(b.1)] If the crossing partition $\crosspart$ is in $P_X$, then $P_X$ is one of the six cases of Theorem \ref{SatzGenauSechs}.
      \item[(b.2)] If neither the crossing partition $\crosspart$ nor the four block partition $\vierpart$ is in $P_X$, then $P_X=\langle \halflibpart\rangle$ or $P_X=\langle\halflibpart,\singleton\otimes\singleton\rangle$.
     \end{itemize}
\end{itemize}
\end{prop}
\Proof (b.2) By Lemma \ref{LemmaD}, all blocks of partitions in $P_X$ are of length at most two. Thus, by Lemma \ref{LemmaC}, the half-liberating partition $\halflibpart$ is in $P_X$. We have $\legpart\;\notin P_X$ by Lemma \ref{LemmaA}.

\emph{Case 1.} If $\singleton\otimes\singleton\;\notin P_X$, then $P_X$ consists of pair partitions. Indeed, if there was a partition in $P_X$ containing a singleton as a block, we could use the pair partition to obtain either $\singleton\in P_X$ or $\singleton\otimes\singleton\in P_X$. Both cases are contradictions. Furthermore, each pair connects one $\oplus$ with one $\ominus$. Otherwise, there would be a pair with an odd number of points between its legs, thus an application of the pair partition $\paarpart$ would yield the crossing partition $\crosspart$ which is a contradiction.
From Proposition \ref{HalfLibKateg}, we infer that $\langle\halflibpart\rangle= P_X$.

\emph{Case 2.} If $\singleton\otimes\singleton\;\in P_X$, then $\langle\halflibpart,\singleton\otimes\singleton\rangle\subset P_X$. On the other hand, every partition in $P_X$ consists of an arbitrary number of blocks of size two and an even number of blocks of size one (a partition of odd length would yield a singleton). Every pair connects one $\oplus$ with one $\ominus$, since otherwise either $\legpart\;\in P_X$ or $\crosspart\;\in P_X$. We conclude $\langle\halflibpart,\singleton\otimes\singleton\rangle= P_X$ using Proposition \ref{HalfLibKateg}.\qqed 

The next theorem is an adaption of \cite[th. 6.5]{BCS10} respecting the two missing easy quantum groups $\Bns$ and $\Bnhl$. An easy quantum group is called \emph{hyperoctahedral}, if its corresponding category of partitions contains the four block partition $\vierpart$ but not the partition $\singleton\otimes\singleton$. Otherwise, it is called \emph{nonhyperoctahedral}.

\begin{thm}\label{CorNonHyp}
 There are exactly 13 nonhyperoctahedral easy quantum groups, namely:
\begin{itemize}
\item $O_n$, $\Onhl$ and $\On$
\item $S_n$ and $\Sn$
\item $B_n$ and $\Bn$
\item $S_n'$ and $\Snp$
\item $B_n'$ and $\Bnp$
\item $\Bnhl$ and $\Bns$
\end{itemize}
\end{thm}

\Proof Let $P_X$ be a nonhyperoctahedral full category of partitions.

\emph{Case 1.} Let $P_X\subset NC$. Then $P_X$ is one of the six nonhyperoctahedral categories of Proposition \ref{PropKategorien} (by Theorem \ref{SatzGenauSieben}).

\emph{Case 2a.} Let $P_X\not\subset NC$ and $\crosspart\;\in P_X$. Then $P_X$ is one of the five nonhyperoctahedral categories of Theorem \ref{SatzGenauSechs}.

\emph{Case 2b.} Let $P_X\not\subset NC$ and $\crosspart\;\notin P_X$. 
Then $\vierpart\;\notin P_X$ by Corollary \ref{CorDoubleSingl}, since $P_X$ is nonhyperoctahedral.
By Proposition \ref{PropNonhyperoct}(b.2), we conclude that $P_X$ is either $\langle \halflibpart\rangle$ or $\langle\halflibpart,\singleton\otimes\singleton\rangle$.\qqed

The remaining easy quantum groups are of hyperoctahedral type and there are many examples of them. It is an open question how to classify them and how to give a complete classification of all easy quantum groups. 
At least, to shed some more light on this, we will now list all categories of \emph{half-liberated type}, i.e. all categories of partitions involving the half-liberating partition $\halflibpart$. For this, we use Proposition \ref{PropNonhyperoct} and the following lemma. 

For $s\in \N$, we denote by $h_s\in P(0,2s)$ the partition given by the two blocks $\{1,3,5,\ldots,2s-1\}$ and $\{2,4,6,\ldots,2s\}$ of length $s$ respectively. 

\setlength{\unitlength}{0.5cm}
\begin{center}
\begin{picture}(12,3)
\put(-1,0){$h_s\;=$}
\put(-0.1,0){\uppartiii{2}{1}{3}{5}}
\put(5,2){\line(1,0){1}}
\put(6.75,1.95){$\ldots$}
\put(-0.1,0){\uppartiii{1}{2}{4}{6}}
\put(6,1){\line(1,0){1}}
\put(7.75,0.95){$\ldots$}
\put(-0.1,0){\uppartii{1}{10}{12}}
\put(8.65,2){\line(1,0){0.5}}
\put(-0.1,0){\uppartii{2}{9}{11}}
\put(9.5,1){\line(1,0){1}}
\end{picture}
\end{center}

For $s=1$, we have $h_1=\;\singleton\otimes\singleton$, for $s=2$ the partition $h_2$ is a rotated version of the crossing partition $\crosspart$.
The partitions $h_s$ were used by Banica, Curran and Speicher in \cite[def. 3.1]{BCS10} to define a \emph{hyperoctahedral series} $H_n^{(s)}$. Their corresponding categories of partitions are given by $\langle\halflibpart, \vierpart, h_s\rangle$, for $s\geq3$. (cf. also \cite[lem. 3.2]{BCS10})

\begin{lem}\label{LemHs}
 Let $P_X$ be a full category of partitions. Assume $\vierpart\;\in P_X$, $\crosspart\;\notin P_X$, $P_X\not\subset NC$ and $P_X\not\subset \langle\halflibpart, \vierpart\rangle$. Then the following statements hold true.
\begin{itemize}
 \item[(a)] There is a $h_s\in P_X$ for some $s\geq 3$.
 \item[(b)] If $h_s\in P_X$ for a $s\geq 3$, then $h_{ks}\in P_X$ for all $k\in \N$.
 \item[(c)] If $h_s\in P_X$ and $h_{s'}\in P_X$ for some $s>s'\geq 3$, then also $h_{s-s'}\in P_X$. 
 \item[(d)] Let $h_s\in P_X$ and $h_{s'}\in P_X$ for some $s,s'\geq 3$. Denote by $g\in \N$ the greatest common divisor of $s$ and $s'$. Then $h_g\in P_X$.
\end{itemize}
\end{lem}
\Proof (a) Since $P_X\not\subset\langle\halflibpart, \vierpart\rangle$, there is a partition $p\in P_X$ containing a block $b\subset p$ which consists of $\alpha$ points labeled by $\oplus$ and $\beta$ points labeled by $\ominus$ such that $\alpha\neq \beta$ (by Proposition \ref{HalfLibKateg}). We may assume $\alpha > \beta$. If there is a  $\oplus$-point and a $\ominus$-point  belonging to $b$ such that all points in between are \emph{not} part of $b$, we erase them (it is an even number) using the pair partition $\paarpart$  before we erase these two points $\oplus$ and $\ominus$, too. 
By this procedure, we may erase all $\ominus$-points of $p$.
We also erase the points in between two $\oplus$-points of $b$, such that only one point is left over respectively.
By iteration, we see that the following partition $p'$ is in $P_X$ where $b'$ denotes the block according to $b$.
\setlength{\unitlength}{0.5cm}
\begin{center}
\begin{picture}(11,3)
\put(-1,1){$p'\;=$}
\put(-0.1,1){\uppartiii{1}{1}{3}{5}}
\put(5,2){\line(1,0){1}}
\put(6.25,1.95){$\ldots$}
\put(-0.1,1){\upparti{1}{10}}
\put(8.5,1.95){$\ldots$}
\put(9.65,2){\line(1,0){0.5}}
\put(2,2.5){$b'$}
\put(0.9,0){$\oplus$}
\put(2.9,0){$\oplus$}
\put(4.9,0){$\oplus$}
\put(9.9,0){$\oplus$}
\put(2,1){$\bullet$}
\put(4,1){$\bullet$}
\put(6,1){$\bullet$}
\put(7.25,1){$\ldots$}
\put(9,1){$\bullet$}
\put(11,1){$\bullet$}
\end{picture}
\end{center}
The block $b'$ is of length at least three. Indeed, $b'$ cannot be a singleton, because $\singleton\otimes\singleton\;\notin P_X$. (by Corollary \ref{CorDoubleSingl}; furthermore note that $p'$ is of even length, since otherwise $\singleton\;\in P_X$.) This also shows that there are no singletons at all in $p'$. Therefore, if $b'$ was of length two, then $p'$ would be a rotated version of the crossing partition $\crosspart$, which is a contradiction to $\crosspart\;\notin P_X$. 

Thus, $b'$ is of length at least three and the points in between the legs of $b'$ are also connected by blocks of length at least three. (If two points were connected by a pair, we would obtain the crossing partition $\crosspart$ in $P_X$  using the pair partition applied to points not belonging to the pair.) Hence, we use the pair partition to obtain $h_s$ in $P_X$, for some $s\geq 3$.

(b) If $h_s\in P_X$ and $h_{ks}\in P_X$, then the partition $h_s\otimes\vierpart\otimes h_{ks}$ is in $P_X$, too. We apply the pair partition $\paarpart$ twice to connect the four block partition $\vierpart$ to $h_s$ and to $h_{ks}$ respectively (i.e. we compose $h_s\otimes\vierpart\otimes h_{ks}$ with $\idpart^{\otimes(2s-1)}\otimes\paarpart^*\otimes\idpart\otimes\idpart\otimes\paarpart^*\otimes\idpart^{\otimes(2ks-1)}$). By this, the following partition $q$ is in $P_X$.
\setlength{\unitlength}{0.5cm}
\begin{center}
\begin{picture}(21,5)
\put(-1,2){$q\;=$}
\put(-0.1,2){\uppartiii{2}{1}{3}{5}}
\put(5,4){\line(1,0){1}}
\put(7.25,3.95){$\ldots$}
\put(-0.1,2){\uppartiii{1}{2}{4}{6}}
\put(6,3){\line(1,0){1}}
\put(7.25,2.95){$\ldots$}
\put(-0.1,2){\uppartiv{1}{9}{11}{12}{14}}
\put(9.65,4){\line(1,0){0.5}}
\put(13.15,4){\line(1,0){0.5}}
\put(15.25,3.95){$\ldots$}
\put(-0.1,2){\upparti{2}{10}}
\put(-0.1,2){\upparti{2}{13}}
\put(8.5,3){\line(1,0){1}}
\put(14,3){\line(1,0){1}}
\put(15.25,2.95){$\ldots$}
\put(-0.1,2){\uppartiii{2}{18}{20}{22}}
\put(17.5,4){\line(1,0){1}}
\put(-0.1,2){\uppartiii{1}{17}{19}{21}}
\put(16.5,3){\line(1,0){1}}
\put(1.1,1.5){$\underbrace{\qquad\qquad\qquad\qquad\qquad\qquad\;\;}$}
\put(5,0){$2s$ points}
\put(12.1,1.5){$\underbrace{\qquad\qquad\qquad\qquad\qquad\qquad\;\;}$}
\put(15,0){$2ks$ points}
\end{picture}
\end{center}
We rotate the first $s$ points of $h_s$ to obtain a partition $r_0\in P(s,s)$. This gives rise to a partition $r:=r_0\otimes r_0\in P(2s,2s)$ which acts like $\crosspart^{\otimes s}$ when applied to $h_s$. Thus, if we compose $q$ with $r^{\otimes (k+1)}$, the following partition $q'$ is in $P_X$.
\setlength{\unitlength}{0.5cm}
\begin{center}
\begin{picture}(21,5)
\put(-1,2){$q'\;=$}
\put(-0.1,2){\uppartiii{2}{2}{4}{6}}
\put(5,3){\line(1,0){2}}
\put(7.25,3.95){$\ldots$}
\put(-0.1,2){\uppartiii{1}{1}{3}{5}}
\put(6,4){\line(1,0){1}}
\put(7.25,2.95){$\ldots$}
\put(-0.1,2){\uppartii{1}{10}{13}}
\put(8.65,4){\line(1,0){0.5}}
\put(14.15,4){\line(1,0){0.5}}
\put(15.25,3.95){$\ldots$}
\put(-0.1,2){\uppartii{2}{9}{11}}
\put(-0.1,2){\uppartii{2}{12}{14}}
\put(8.65,3){\line(1,0){2}}
\put(13.15,3){\line(1,0){1.5}}
\put(15.25,2.95){$\ldots$}
\put(-0.1,2){\uppartiii{2}{17}{19}{21}}
\put(16.5,3){\line(1,0){2}}
\put(-0.1,2){\uppartiii{1}{18}{20}{22}}
\put(16.5,4){\line(1,0){1}}
\put(1.1,1.5){$\underbrace{\qquad\qquad\qquad\qquad\qquad\qquad\;\;}$}
\put(5,0){$2s$ points}
\put(12.1,1.5){$\underbrace{\qquad\qquad\qquad\qquad\qquad\qquad\;\;}$}
\put(15,0){$2ks$ points}
\end{picture}
\end{center}
We compose $q'$ with $\idpart^{\otimes 2s}\otimes\vierpart\otimes\idpart^{\otimes 2ks}$ and use the pair partition twice to connect the four block partition $\vierpart$ to $q'$. Thus, the resulting partition $q''$ is in $P_X$.
\setlength{\unitlength}{0.5cm}
\begin{center}
\begin{picture}(21,5)
\put(-1,2){$q''\;=$}
\put(-0.1,2){\uppartiii{2}{2}{4}{6}}
\put(5,3){\line(1,0){2}}
\put(7.25,3.95){$\ldots$}
\put(-0.1,2){\uppartiii{1}{1}{3}{5}}
\put(6,4){\line(1,0){1}}
\put(7.25,2.95){$\ldots$}
\put(-0.1,2){\uppartii{1}{10}{13}}
\put(8.65,4){\line(1,0){0.5}}
\put(14.15,4){\line(1,0){0.5}}
\put(15.25,3.95){$\ldots$}
\put(-0.1,2){\uppartiv{2}{9}{11}{12}{14}}
\put(8.65,3){\line(1,0){2}}
\put(13.15,3){\line(1,0){1.5}}
\put(15.25,2.95){$\ldots$}
\put(-0.1,2){\uppartiii{2}{17}{19}{21}}
\put(16.5,3){\line(1,0){2}}
\put(-0.1,2){\uppartiii{1}{18}{20}{22}}
\put(16.5,4){\line(1,0){1}}
\put(1.1,1.5){$\underbrace{\qquad\qquad\qquad\qquad\qquad\qquad\;\;}$}
\put(5,0){$2s$ points}
\put(12.1,1.5){$\underbrace{\qquad\qquad\qquad\qquad\qquad\qquad\;\;}$}
\put(15,0){$2ks$ points}
\end{picture}
\end{center}
An application of the partition $r$ to the first $2s$ points of $q''$ yields the partition $h_{(k+1)s}$.

(c) Compose $h_s$ with $h_{s'}^*\otimes \idpart^{\otimes 2(s-s')}$ to obtain $h_{s-s'}$.

(d) We deduce (d) from (b) and (c) by some arithmetics. Denote by $B_X\subset \N$ the set of all $t\geq 3$, such that $h_t\in P_X$. Thus, $B_X$ has the following properties:
\begin{itemize}
 \item The numbers $s$ and $s'$ are in $B_X$.
\item If $t\in B_X$ and $k\in \N$, then $kt\in B_X$.
\item If $t,t'\in B_X$ such that  $t>t'$, then $t-t'\in B_X$.
\item If $ag, bg\in B_X$ such that $a>b$ and $a\notin b\N$, then there is a $b'<b$ such that $b'g\in B_X$.
\end{itemize}
To prove the last item, take $k\in \N$ such that $bk<a<b(k+1)$ and put $b':=a-bk$.

Now, we write $s=ag$ and $s'=bg$ where $g$ is the greatest common divisor. We may assume $a>b$, thus there is a $b'<b$ such that $b'g\in B_X$. Then $b\notin b'\N$ or $a\notin b'\N$ or $b'=1$. If $b'\neq 1$, then again there is a $b''<b'$ such that $b''g\in B_X$ and $b\notin b''\N$ or $a\notin b''\N$ or $b''=1$. By induction, we conclude $g\in B_X$.
\qqed

Next, we list all categories of half-liberated type. 

\begin{thm}
 Let $P_X$ be a full category of partitions containing the half-liberating partition $\halflibpart$, but not the crossing partition $\crosspart$. Then $P_X$ is one of the following cases:
\begin{itemize}
 \item $P_X=\langle\halflibpart\rangle$, corresponding to $\Onhl$.
 \item $P_X=\langle\halflibpart, \singleton\otimes\singleton\rangle$, corresponding to $\Bnhl$.
\item $P_X=\langle\halflibpart, \vierpart\rangle$, corresponding to $\Hnhl$.
 \item $P_X=\langle\halflibpart, \vierpart, h_s\rangle$ for some $s\geq 3$, corresponding to the hyperoctahedral series $H_n^{(s)}$.
\end{itemize}
\end{thm}
\Proof By Proposition \ref{PropNonhyperoct}(b.2), we can restrict to the case $\vierpart\;\in P_X$. If \linebreak $P_X\neq\langle\halflibpart, \vierpart\rangle$, there is a $h_t\in P_X$ by Lemma \ref{LemHs} for some $t\geq 3$. Denote by $s$ the minimal number $s\geq 3$ such that $h_s\in P_X$. Then $P_X\supset\langle\halflibpart, \vierpart, h_s\rangle$. For the converse direction, let $p\in P_X$ and let $b\subset p$ be a block of $p$. Since the half-liberating partition $\halflibpart$ is in $P_X$, we can arrange the legs of $b$ such that $p$ is of the following form. (We may assume that the number of points of $b$ labeled by $\oplus$ is greater or equal the number of points labeled by $\ominus$.)
\setlength{\unitlength}{0.5cm}
\begin{center}
\begin{picture}(18,3)
\put(2,2.5){$b$}
\put(-1,1){$p\;=$}
\put(-0.1,1){\uppartiv{1}{1}{2}{3}{4}}
\put(4,2){\line(1,0){1}}
\put(5.25,1.95){$\ldots$}
\put(0.9,0){$\oplus$}
\put(1.9,0){$\ominus$}
\put(2.9,0){$\oplus$}
\put(3.9,0){$\ominus$}
\put(5.25,0){$\ldots$}
\put(6.5,2){\line(1,0){1}}
\put(-0.1,1){\uppartiv{1}{7}{8}{9}{11}}
\put(11,2){\line(1,0){1}}
\put(10,0.8){$\bullet$}
\put(12,0.8){$\bullet$}
\put(6.9,0){$\oplus$}
\put(7.9,0){$\ominus$}
\put(8.9,0){$\oplus$}
\put(10.9,0){$\oplus$}
\put(12.8,0){$\ldots$}
\put(12.8,2){$\ldots$}
\put(14.5,2){\line(1,0){1}}
\put(-0.1,1){\uppartii{1}{15}{17}}
\put(14,0.8){$\bullet$}
\put(16,0.8){$\bullet$}
\put(18,0.8){$\bullet$}
\put(14.9,0){$\oplus$}
\put(16.9,0){$\oplus$}
\put(18.8,0){$\ldots$}
\end{picture}
\end{center}
Here, the bullet points are connected by some other partitions. We use the pair partition $\paarpart$ to erase all $\ominus$-points of $b$ and we obtain a partition $p'\in P_X$. Then $p\in \langle\halflibpart, \vierpart, h_s\rangle$ if and only if $p'\in \langle\halflibpart, \vierpart, h_s\rangle$, since we may apply the pair partition to $\vierpart^{\otimes k}\otimes p'$ to reconstruct $p$. Iterating this procedure yields a partition $q\in P_X$, whose blocks are labeled either only by $\oplus$ or only by $\ominus$ respectively. Furthermore, $p\in \langle\halflibpart, \vierpart, h_s\rangle$ if and only if $q\in \langle\halflibpart, \vierpart, h_s\rangle$.

If $t$ is the length of a block $b$ of $q$, then $h_t\in P_X$. Indeed, let $k\in \N$ be large enough such that $ks\geq t$ and assume that all points of $b$ are labeled by $\oplus$. Then $q\otimes h_{ks}$ is in $P_X$ containing a block $r$ of length $ks$ labeled by $\ominus$-points. Using the half-liberating partition $\halflibpart$, we can shift $r$ to the block $b$ in such a way that all $\ominus$-points in between the legs of $b$ are connected by the block $r$. Using the pair partition $\paarpart$ yields that $h_t$ is in $P_X$.

We conclude that the lengths of all blocks of $q$ are multiples of $s$, by Lemma \ref{LemHs} and minimality of $s$. Thus, $q$ is in $\langle\halflibpart, \vierpart, h_s\rangle$ since it may be constructed out of $h_s^{\otimes k}$ for some $k$, using the four block partition $\vierpart$ and the pair partition $\paarpart$ to connect some blocks (see also the proof of Lemma \ref{LemHs}). Hence, $p\in \langle\halflibpart, \vierpart, h_s\rangle$.\qqed

\begin{cor}\label{CorHnHnhl}
 If $G$ is an easy quantum group with $H_n\subsetneqq G\subsetneqq\Hnhl$, then $G=H_n^{(s)}$ for some $s\geq 3$.
\end{cor}


By Theorem \ref{CorNonHyp}, we know that there are exactly 13 nonhyperoctahedral easy quantum groups. The remaining easy quantum groups are all of hyperoctahedral type, i.e. their categories of partitions contain the four block partition $\vierpart$ but not the double singleton $\singleton\otimes\singleton$.
The complete classification in the hyperoctahedral case has to be left open. 
Nevertheless, we briefly sketch the state of the art introducing a further example of a (hyperoctahedral) easy quantum group.

There are two more types of partitions that need to be considered.
For a natural number $l\in\N$, we denote by $k_l\in P(l+2,l+2)$ the partition given by a four block on $\{1,1',l+2,(l+2)'\}$ and pairs on $\{i,i'\}$ for $i=2,\ldots,l+1$. (cf. also \cite[lem. 4.2]{BCS10}) The following picture illustrates the partition $k_l$ -- note that the waved line from $1'$ to $l+2$ is \emph{not} connected to the lines from $2$ to $2'$ from $3$ to $3'$ etc.
\setlength{\unitlength}{0.5cm}
\begin{center}
\begin{picture}(14,5)
\put(0.5,2.5){$k_l\;\;=$}
\put(1.9,1.5){\upparti{2}{1}}
\put(1.9,1.5){\upparti{2}{2}}
\put(1.9,1.5){\upparti{2}{3}}
\put(7,2){$\ldots$}
\put(1.9,1.5){\upparti{2}{8}}
\put(1.9,1.5){\upparti{2}{9}}
\put(3,0){$1'$}
\put(4,0){$2'$}
\put(5,0){$3'$}
\put(7,0){$\ldots$}
\put(8.4,0){$(l+1)'$}
\put(11,0){$(l+2)'$}
\put(3,4){$1$}
\put(4,4){$2$}
\put(5,4){$3$}
\put(7,4){$\ldots$}
\put(9.1,4){$l+1$}
\put(10.9,4){$l+2$}
\put(7.2,1.5){\oval(8,2)[tl]}
\put(7.15,3.5){\oval(8,2)[br]}
\end{picture}
\end{center}

Secondly, we consider the \emph{fat crossing partition} $\fatcrosspart\;\in P(4,4)$ given by $\{1,2,3',4'\}$ and $\{3,4,1',2'\}$. It is not just a double crossing but two crossing four block partitions.
In rotated form, it is of the following form:
\setlength{\unitlength}{0.5cm}
\begin{center}
\begin{picture}(8,3)
\put(0.5,0){$\fatcrosspart\;\;=$}
\put(1,0){\uppartiv{1}{1}{2}{5}{6}}
\put(1,0){\uppartiv{2}{3}{4}{7}{8}}
\end{picture}
\end{center}


We have $k_1\in\langle\halflibpart,\vierpart\rangle$ as well as $\fatcrosspart\;\in\langle k_1,\vierpart\rangle$ and $k_1, \fatcrosspart\;\in\langle \vierpart,h_s\rangle$ for all $s\geq 3$. Furthermore, a category contains the partitions $k_l$ for all $l\in\N$ if and only if it contains $k_1$ (apply the pair partition to $k_l\otimes k_1$ to obtain $k_{l+1}$).
Currently, we have the following list of hyperoctahedral categories of partitions, including a new one.
\begin{itemize}
 \item $\langle\vierpart\rangle$, $\langle\halflibpart,\vierpart\rangle$ and $\langle\crosspart,\vierpart\rangle$ corresponding to $\Hn$, $\Hnhl$ and $H_n$ respectively.
 \item $\langle \vierpart,h_s\rangle$ for $3\leq s<\infty$ and $\langle \vierpart, k_1\rangle$ for $s=\infty$ corresponding to the \emph{higher hyperoctahedral series} $\Hnseckig$, defined in \cite{BCS10}.
 \item The series $\langle\halflibpart,\vierpart,h_s\rangle$ parametrized by $s\geq 3$ corresponding to $\Hnsrund$.
 \item The new category $\langle\fatcrosspart,\vierpart\rangle$.
\end{itemize}

The partition $\fatcrosspart$ gives rise to the relations $u_{ij}^2u_{kl}^2=u_{kl}^2u_{ij}^2$.
Hence, for $n\in \N$, the new quantum group is given by:
\[\univ{u_{ij}, i,j=1,\ldots,n}{(u_{ij}) \textnormal{ is cubic and } u_{ij}^2u_{kl}^2=u_{kl}^2u_{ij}^2}\]
We will study the hyperoctahedral case of easy quantum groups in a future article.

\section{Laws of characters}\label{Sect5}

An important step in the investigation of an easy quantum group $G_n$ with associated category of partitions $P_{G_n}$ is the computation of the law of the character $\chi=\sum_{i=1}^n u_{ii}\in G_n$. (cf. \cite[sect. 5]{BS09}) The question is to find a real probability measure $\mu$ such that for all $k\in \N$:
\[\int_\R x^k\textnormal{d}\mu(x)=\int \chi ^k\]
Here, $\int$ denotes the Haar functional on the quantum group by Woronowicz (\cite{Wo87}). 
For easy quantum groups, we can solve this problem in an asymptotic way using the following formula from Banica, Curran and Speicher (\cite[th. 2.3]{BCS11}):
\[\lim_{n\to\infty} \int\chi^k=\#\{p\in P(0,k)\;|\;p\in P_{G_n}\}=:m_k\]
Note, that if $P_X$ is a category of partitions, in fact we assign a \emph{family} of easy quantum groups $G_n$ to it such that $P_{G_n}=P_X$, for all dimensions $n\in\N$. For a fixed  $n\in\N$, the law of characters may be computed only up to lower order correction terms, using a Weingarten formula (see \cite{BCS11}). In the limit, we get the above exact formula in terms of the category $P_X=P_{G_n}$, representing the asymptotic law of characters.

We may also express the asymptotic law of $\chi\in G_n$ in terms of elements $a\in B$ in a  *-probability space $(B,\phi)$, if $\phi(a^k)=m_k$ for all $k$.
We say that the law of $\chi$ is given by the \emph{squeezed version} of $a$ (denoted by $\sqrt{aa^*}$), if we are in the following situation. (cf. \cite[def. 7.5]{BCS10}) 
\[m_{l}=\begin{cases}\phi((aa^*)^k) \quad&\textnormal{if } l=2k\\ 0 \quad &\textnormal{otherwise}\end{cases}\]
 
Recall some basic facts from free probability theory (see for instance \cite{NiSp}). A \emph{semicircular element} $s$ is characterized by its free cumulants:
\[\kappa_2(s,s)=1,\qquad \kappa_n(s,\ldots,s)=0, \;n\neq 2\]
A \emph{circular element} $c$ is characterized by:
\[\kappa_2(c,c^*)=\kappa(c^*,c)=1,\qquad \textnormal{all other cumulants vanish}\]
Circular elements are complex versions of semicircular elements, since\linebreak $c=\frac{1}{\sqrt 2}(s_1+is_2)$ is circular, if $s_1$ and $s_2$ are free semicircular elements. The semicircular elements are free versions of the real Gaussian law $g$, whereas circular elements correspond to the complex Gaussian law $\tilde g$.
This correspondence is also reflected by the laws of $O_n$, $\Onhl$ and $\On$.

\begin{prop}[{\cite{BBC07},\cite{BS09},\cite{BCS10}}] 
The asymptotic laws of the classical, half-liberated and free versions of the orthogonal group are:
\begin{itemize}
 \item The law of $O_n$ is the real Gaussian.
\item The law of $\Onhl$ is the squeezed complex Gaussian.
\item The law of $\On$ is the semicircle.
\end{itemize}
\end{prop}
Note that the squeezed version of a circular element is a semicircular element. We illustrate the above correspondence by the following table.

\begin{center}\begin{tabular}{rclrc}
$O_n$: &$g$ &\qquad\qquad\qquad
&$\On$: &$s$ \\
        &real Gaussian&
&       &semicircular\\
\rule[0mm]{0mm}{6mm}
$\Onhl$: &$\sqrt{\tilde g\tilde g^*}$ \\
        &squeezed\\
        &complex Gaussian
\end{tabular}
\end{center}

In the bistochastic' case, we have a similar correspondence which explains the splitting of $B_n'$ into two free versions (cf. Remark \ref{RemY}) from the free probability point of view. The \emph{shifted} version of a law $x$ is the law of $1+x$. The \emph{symmetrized} version $\mu'$ of a measure $\mu$ is given by $\mu'(B)=\frac{1}{2}(\mu(B)+\mu(-B))$. It may also be expressed in terms of moments. The symmetrized version of a law $x$ is the law $x'$ with the same even moments and vanishing odd moments. Thus, the squeezed version of a self-adjoint element coincides with its symmetrized version.

\begin{prop}\label{PropBistochpLaws} 
The asymptotic laws of the several variants of the bistochastic group are:
\begin{itemize}
 \item The law of $B_n$ is a shifted real Gaussian.
 \item The law of $B_n'$ is a symmetrized shifted real Gaussian (which is the squeezed version of a shifted real Gaussian).
 \item The law of $\Bnhl$ is the squeezed version of a shifted complex Gaussian.
 \item The law of $\Bn$ is a shifted semicircle.
 \item The law of $\Bnp$ is a symmetrized shifted semicircle (which is the squeezed version of a shifted semicircle).
 \item The law of $\Bns$ is  the squeezed version of a shifted circle.
\end{itemize}
 \end{prop}
\Proof The laws of $B_n$, $B_n'$, $\Bn$ and $\Bnp$ have been computed in \cite[th. 5.6]{BS09}.

For the computation of the law of $\Bns$, we put $d:=1+c$ and we check:
\[ \kappa_2(d,d^*)=\kappa_2(d^*,d)=\kappa_1(d)=\kappa_1(d^*)=1\]
All other cumulants of $d$ vanish. By the moment-cumulant formula (see for instance \cite{NiSp}), the $k$-th moment of $dd^*$ is of the form:
\[\phi((dd^*)^k)=\sum_{\pi\in NC(0,2k)} \kappa_\pi(d,d^*,\ldots,d,d^*)\]
Now, $\kappa_\pi(d,d^*,\ldots,d,d^*)=1$ if and only if $\pi\in\langle\singleton\otimes\singleton\rangle$ and zero otherwise (using Proposition \ref{PropKategorien}). Hence, $\phi((dd^*)^k)$ is exactly the number of partitions in $\langle\singleton\otimes\singleton\rangle$ of length $2k$.

The element $\tilde d:=1+\tilde g$ has the same (classical) cumulants as $d$ and by the moment-cumulant formula for classical cumulants, we infer the statement for $\Bnhl$. (Compare the categories of $\Bns$ in Proposition \ref{PropKategorien} and $\Bnhl$ in Proposition \ref{HalfLibKateg}.)
 \qqed

The preceding proposition is illustrated by the following table.
\begin{center}\begin{tabular}{rclrc}
$B_n$: &$1+g$ &\qquad\qquad\qquad
&$\Bn$: &$1+s$ \\
        &shifted &
&       &shifted\\
        &real Gaussian&
&       &semicircular\\
\rule[0mm]{0mm}{6mm}
$B_n'$: &$\sqrt{(1+g)(1+g)^*}$ &
&$\Bnp$: &$\sqrt{(1+s)(1+s)^*}$ \\
        &squeezed shifted &
&       &squeezed shifted\\
        &real Gaussian &
&       &semicircular\\
\rule[0mm]{0mm}{6mm}
$\Bnhl$: &$\sqrt{(1+\tilde g)(1+\tilde g)^*}$ &
&$\Bns$: &$\sqrt{(1+c)(1+c)^*}$ \\
        &squeezed shifted &
&       &squeezed shifted\\
        &complex Gaussian &
&       &circular
\end{tabular}
\end{center}

Thus, the laws of the quantum groups $\Bnhl$ and $\Bns$ appear naturally as ``complex'' versions of the laws of $B_n'$ and $\Bnp$. The analogy to the table of the non-shifted laws is evident -- since the real Gaussian and the semicircular are already symmetric, there are no symmetrized versions of them.

We end this section by the computation of the moments of $\Bns$. It is due to Octavio Arizmendi (private communication). 

\begin{prop}\label{RemOctavio}
 The odd moments of the asymptotic law of the character $\chi=\sum_{i=1}^nu_{ii}$ in $\Bns$ are all zero and the even moments of order $2k$ are given by
\[b_k=\frac{1}{k+1}\binom{3k+1}{k}\]
The sequence $(b_k)_{k\in \N}$ is related to the (type 2) Fuss-Catalan numbers given by
\[C_k^{(2)}=\frac{1}{2k+1}\binom{3k}{k}\]
Let $g(x):=\sum_{k=0}^\infty C_k^{(2)}x^k$, then $g^2(x)=\sum_{k=0}^\infty b_kx^k$.
(cf. \cite{GX05})
\end{prop} 
\Proof We have to count the number of partitions $p\in NC(0,m)$ in $\langle\singleton\otimes\singleton\rangle$ of length $m$. 
Take a partition $p\in\langle\singleton\otimes\singleton\rangle$ of length $2k$. We compose a pair partition $\paarpart$ with $\idpart\otimes p\otimes\idpart$ to obtain a partition $p'\in NC(0,2k+2)$. Now, we connect all singletons of $p'$ with the next pair partition ``above'' it (this is the reason, why we consider $p'$ instead of $p$ itself). This may be illustrated by the following picture.
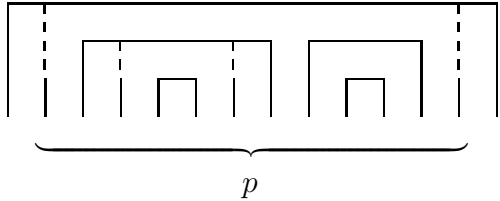
\setlength{\unitlength}{0.5cm}
\begin{center}
\begin{picture}(15,6)
\put(0,2){\uppartii{3}{1}{14}}
\put(0,2){\upparti{1}{2}}
\put(2.25,3){\usebox{\dashlinehoch}}
\put(2.25,4){\usebox{\dashlinehoch}}
\put(0,2){\uppartii{2}{3}{8}}
\put(0,2){\upparti{1}{4}}
\put(4.25,3){\usebox{\dashlinehoch}}
\put(0,2){\uppartii{1}{5}{6}}
\put(0,2){\upparti{1}{7}}
\put(7.25,3){\usebox{\dashlinehoch}}
\put(0,2){\uppartii{2}{9}{12}}
\put(0,2){\uppartii{1}{10}{11}}
\put(0,2){\upparti{1}{13}}
\put(13.25,3){\usebox{\dashlinehoch}}
\put(13.25,4){\usebox{\dashlinehoch}}
\put(2,1.5){$\underbrace{\qquad\qquad\qquad\qquad\qquad\qquad\qquad}$}
\put(7.5,0){$p$}
\end{picture}
\end{center}
Hence, the partitions in $\langle\singleton\otimes\singleton\rangle$ of length $2k$ are in one-to-one correspondence to the even partitions in $NC(0,2k+2)$, where $1$ and $2k+2$ are in the same block. These in turn are in one-to-one correspondence to the 3-equal partitions of length $3k+3$ where 1 and $3k+3$ are in the same block. This number $b_k$ may be determined explicitly.
(cf. also \cite[ex. 22]{A}) \qqed

\section{About the structure of the $C^*$-algebras associated to the free easy quantum groups}\label{Sect3}

In order to understand the difference between the two quantum versions $\Bnp$ and $\Bns$ of $B_n'$, it is useful to consider only the $C^*$-algebraic structure of these objects. First, we will prove that the $C^*$-algebras $\ABns$ and $\ABnp$ may be constructed out of $\ABn$ in a direct way, namely by the tensor product with $C^*(\Z_2)$ respectively the free product. Analogously, $\ASnp$ may be obtained out of $\ASn$. Secondly, we will deduce some consequences for the $K$-theory of these $C^*$-algebras, before we investigate the property of exactness.

Let us recall some $C^*$-algebraic constructions.
We view the tensor product $A\otimes B$ of a unital, nuclear $C^*$-algebra $A$ with an arbitrary unital $C^*$-algebra $B$ as the universal $C^*$-algebra generated by all elements $a\in A$ (with the relations of $A$) and all $b\in B$ (with the relations of $B$), such that all $a\in A$ commute with all $b\in B$. The units of $A$ and $B$ are identified. (cf. also \cite{Cu93})

Likewise, the unital free product $A *_\C B$ of two unital $C^*$-algebras $A$ and $B$ is the universal $C^*$-algebra generated by all elements $a\in A$ (again with the relations of $A$) and all $b\in B$ (with the relations of $B$) without any further relations. The units of $A$ and $B$ are identified.

By $\Z_2$ we denote the group of integers $\Z$ modulo $2\Z$. Its full group $C^*$-algebra $C^*(\Z_2)$ may be seen as the universal $C^*$-algebra generated by a symmetry $s$, i.e. $s^*=s$ and $s^2=1$. Indeed, the generator $a$ of the group $\Z_2$ fulfills $a^2=1$, which caries over to the generating unitary $s$ of $C^*(\Z_2)$. Thus, $C^*(\Z_2)$ is a two-dimensional commutative $C^*$-algebra, and hence it is isomorphic to $\C^2$.

The idea of the isomorphism of $\ABnp$ and $\ABn\otimes C^*(\Z_2)$ is to extract the  symmetry $r=\sum_{k=1}^n u_{ik}\in\ABnp$ (cf. Remark \ref{BemUeberR}).
The product of $r$ with the $u_{ij}$ acts like a tensor product since $r$ commutes with all $u_{ij}$. Thus, we map $r$ to the symmetry $s\in C^*(\Z_2)$. The idea of $\ASnp\cong\ASn\otimes C^*(\Z_2)$ is the same.

We will use the same notation, namely $u_{ij}$, for the generators in $\ABnp$, $\ABn$, $\ASnp$ and $\ASn$ respectively. Hopefully this will not cause any confusions.

\begin{prop}\label{PropKonstrIsomI}
 The $C^*$-algebras $\ABnp$ and $\ASnp$ may be constructed as follows, if $n\geq 2$.
\begin{itemize}
 \item[(a)] $\ABnp$ is isomorphic to $\ABn\otimes C^*(\Z_2)$ via $u_{ij}\mapsto u_{ij}s$.
 \item[(b)] $\ASnp$ is isomorphic to $\ASn\otimes C^*(\Z_2)$ via $u_{ij}\mapsto u_{ij}s$.
\end{itemize}
\end{prop}
\Proof (a) The $C^*$-algebra $\ABn\otimes C^*(\Z_2)$ is the universal $C^*$-algebra generated by self-adjoint elements $u_{ij}$, for $i,j=1,\ldots,n$, which fulfill the bistochastic relations of Definition \ref{DefSieben}, together with a symmetry $s$, such that all $u_{ij}$ commute with $s$. The elements $u_{ij}':=u_{ij}s\in \ABn\otimes C^*(\Z_2)$ are self-adjoint and they satisfy the bistochastic' relations. Indeed, check for instance that $\sum_{k=1}^n u_{ik}'u_{jk}'=\sum_{k=1}^n u_{ik}u_{jk}=\delta_{ij}$ and check furthermore $\sum_{k=1}^n u_{ik}'=s$, which commutes with the $u_{ij}'$.
Hence, there is a homomorphism mapping $u_{ij}\mapsto u_{ij}'$ by the universal property. 

For the converse direction, verify that the $u_{ij}'':=u_{ij}r$ and $s'':=r$ in $\ABnp$ satisfy the relations of $\ABn\otimes C^*(\Z_2)$. Recall that $r$ denotes the sum $\sum_{k=1}^n u_{ik}\in\ABnp$, commuting with all $u_{ij}$, and that it is a symmetry (cf. Remark \ref{BemUeberR}). Check that $\sum_{k=1}^n u_{ik}''=r^2=1$.
Thus, we get a homomorphism which is inverse to the above one. 

(b) In the same way, we get the statement for (b). This is, we check again that the elements $u_{ij}':=u_{ij}s\in \Sn\otimes C^*(\Z_2)$ fulfill the orthogonality relations and $\sum_{k=1}^n u_{ik}'=s$. Furthermore $u_{ik}'u_{jk}'=u_{ik}u_{jk}=0$ if $i\neq j$. For the converse direction, note that the elements $u_{ij}'':=u_{ij}r$ and $s'':=r$ fulfill $\sum_{k=1}^n u_{ik}''=r^2=1$.
 \qqed

Proposition \ref{PropKonstrIsomI} may be read as $\Bnp\simeq\Bn\times\Z_2$ and $\Snp\simeq\Sn\times\Z_2$ in analogy to the classical situation of groups (cf. \cite[prop. 2.4]{BS09}).

In $\ABns$, the symmetry $r$ does \emph{not} commute with the $u_{ij}$ and their product behaves more like a free product. Also, we cannot map $u_{ij}$ of $\ABns$ to $u_{ij}s$ of $\ABn *_\C C^*(\Z_2)$, since $u_{ij}s$ is not self-adjoint. We have to chose another map, although the idea of mapping $r\in\ABns$ to $s\in C^*(\Z_2)$ is still the same. It turns out that the isomorphism of $\ABns$ and $\ABn *_\C C^*(\Z_2)$ gives rise to a second isomorphism of $\ABnp$ and \linebreak $\ABn\otimes C^*(\Z_2)$. The next proposition may be read as $\Bns\simeq\Bn * \Z_2$ in some sense.

\begin{prop}\label{PropKonstrIsomII}
 The $C^*$-algebras $\ABnp$ and $\ABns$ may be constructed as follows, if $n\geq 2$.
\begin{itemize}
 \item[(a)] $\ABns$ is isomorphic to $\ABn *_\C C^*(\Z_2)$ via $u_{ij}\mapsto u_{ij}+\frac{1}{n}(s-1)$.
 \item[(b)] $\ABnp$ is isomorphic to $\ABn\otimes C^*(\Z_2)$ via $u_{ij}\mapsto u_{ij}+\frac{1}{n}(s-1)$.
\end{itemize}
\end{prop}
\Proof (a) The elements $u_{ij}':= u_{ij}+\frac{1}{n}(s-1)\in\ABn *_\C C^*(\Z_2)$ are self-adjoint. They also fulfill the orthogonality relations:
\[\sum_{k=1}^n u_{ik}'u_{jk}'=\sum_{k=1}^n u_{ik}u_{jk}  +\frac{1}{n}(s-1)\sum_{k=1}^n u_{jk} +\frac{1}{n}\sum_{k=1}^n u_{ik}(s-1) + \frac{1}{n}(s-1)^2\]
By $\sum_{k=1}^n u_{ik}=\sum_{k=1}^n u_{jk}=1$ and $(s-1)^2=2(1-s)$, we get $\sum_{k=1}^n u_{ik}'u_{jk}'=\sum_{k=1}^n u_{ik}u_{jk}=\delta_{ij}$.  Finally $\sum_{k=1}^n u_{ik}'=s$, thus we obtain a homomorphism from $\ABns$ to $\ABn *_\C C^*(\Z_2)$, mapping $u_{ij}\mapsto u_{ij}'$. On the other hand, a similar computation shows the existence of a homomorphism from $\ABn *_\C C^*(\Z_2)$ to $\ABns$ mapping $u_{ij}$ to $u_{ij} +\frac{1}{n}(1-r)$ and $s$ to $r$. The two homomorphisms are inverse to each other.

(b) Additionally to the computations in (a), we check that the elements \linebreak $u_{ij}':= u_{ij}+\frac{1}{n}(s-1)\in\ABn \otimes C^*(\Z_2)$ commute with the sums $\sum_{k=1}^nu_{ik}'=s$, and the $u_{ij} +\frac{1}{n}(1-r)\in \ABnp$ commute with $r$.  \qqed

The preceding propositions show that $\ABnp$ can be seen as the $C^*$-algebra \linebreak $\ABn\otimes C^*(\Z_2)$ in two ways. This is \emph{not} true for $\ASnp$ and  $\ASn\otimes C^*(\Z_2)$. The isomorphism of Proposition \ref{PropKonstrIsomII} does not carry over to the case of $\ASnp$ and $\ASn\otimes C^*(\Z_2)$. This is, because the elements $u_{ij}':=u_{ij}+\frac{1}{n}(s-1)\in \ASn\otimes C^*(\Z_2)$ do not fulfill $u_{ik}'u_{jk}'=0$ for $i\neq j$. 

Also, in analogy to Proposition \ref{PropKonstrIsomII}, we could consider the $C^*$-algebra \linebreak $\ASn *_\C C^*(\Z_2)$ together with generators $u_{ij}':= u_{ij} + \frac{1}{n}(s-1)\in \ASn *_\C C^*(\Z_2)$. But this does not give rise to an orthogonal quantum group in the sense of Definition \ref{DefOrthQG}, since the map $\Delta(u_{ij}')=\sum_{k=1}^n u_{ik}'\otimes u_{kj}'$ fails to be a homomorphism.

\begin{rem}\label{RemY}
 The two Propositions \ref{PropKonstrIsomI} and \ref{PropKonstrIsomII} explain that $\ABns$ carries much more free structure than $\ABnp$. Also, it shows that $\ASnp$ behaves more like $\ABnp$ rather than like $\ABns$. This justifies why we renamed the quantum group $\Bnp$ of \cite{BS09} to $\Bns$ and why we again used the name $\Bnp$, now for our new quantum group.

In the commutative case (i.e. if the generators $u_{ij}$ commute), the $C^*$-algebras $\ABnp$ and $\ABns$ coincide. Therefore, we only have six easy groups but seven free easy quantum groups. In other words, passing from the easy  groups to the free easy quantum groups, the group $B_n'$ splits into two cases, namely into $\Bnp$ and $\Bns$. On the level of $C^*$-algebras, this difference is explained by the two preceding propositions, i.e. $B_n'$ splits into a ``tensor product version'' and a ``free product version''. See also Proposition \ref{PropBistochpLaws}.
\end{rem}

Independently, Sven Raum (\cite[th. 4.1]{R11}) discovered the statements of Proposition \ref{PropKonstrIsomI}(b) and Proposition \ref{PropKonstrIsomII}(a), although our proofs are different for the latter case. (Besides, note that Raum uses the former definition of $\ABnp$ which is now called $\ABns$.) 

Also, he states that $A_{O_{n-1}^+}$ is isomorphic to $\ABn$, following Banica's and Speicher's proof for the isomorphism of the groups $O_{n-1}$ and $B_n$ (\cite[prop. 2.4]{BS09}). The idea is as follows. Let $A$ be a unital $C^*$-algebra. An orthogonal matrix $u\in M_n(A)$ is bistochastic if and only if $u\xi=\xi$ (the rows sum up to 1) and $u^t\xi=\xi$ (the columns sum up to 1) for the vector $\xi$ filled with 1's. Now, let $T\in M_n(\R)$ be an orthogonal matrix (i.e. $TT^t=T^tT=1$, where $T^t$ is the transposed matrix), such that $Te_1=\frac{1}{\sqrt{n}}\xi$, where $e_1$ is the canonical vector $e_1=(1,0,\ldots,0)$. If $u'\in M_{n-1}(A)$ is orthogonal, then $T\begin{pmatrix} 1 & \\ & u'\end{pmatrix}T^t$ is bistochastic by the above characterization. On the other hand, if $u\in M_n(A)$ is bistochastic, then $T^tuT=\begin{pmatrix} 1 & \\ & u'\end{pmatrix}$, where $u'\in M_{n-1}(A)$ is orthogonal. By this, we get Raum's result.

\begin{prop}[\cite{R11}]\label{PropRaum}
 The $C^*$-algebras $A_{O_{n-1}^+}$ and $\ABn$ are isomorphic, if $n\geq 3$. More precisely, the homomorphisms 
\begin{align*}
&\alpha:A_{O_{n-1}^+}\to\ABn,\quad &u_{ij}\mapsto v_{i+1,j+1}:=\sum_{k,l=1}^n t_{l,i+1}t_{k,j+1}u_{lk}\quad &i,j=1,\ldots,n-1\\
&\beta:\ABn\to A_{O_{n-1}^+},\quad &u_{ij}\mapsto w_{ij}:=\frac{1}{n}+\sum_{k,l=2}^n t_{il}t_{jk}u_{l-1,k-1}\quad &i,j=1,\ldots,n
\end{align*}
are inverse to each other. Here $T=(t_{ij})_{i,j=1,\ldots,n}\in M_n(\R)$ is given as above.
\end{prop}
\Proof The elements $v_{ij}$ resp. $w_{ij}$ are given by the idea described above. For a direct computation, use that $t_{k1}=\frac{1}{\sqrt{n}}$, $k=1,\ldots,n$, and $\sum_{k=1}^n t_{ki}=\sqrt{n}\sum_{k=1}^n t_{ki}t_{k1}=0$ for all $i\neq 1$.\qqed

These isomorphism of $C^*$-algebras enable us to transfer a result of Voigt (\cite{V11}) on the $K$-theory of $\On$ to the quantum groups $\Bn$, $\Bnp$ and $\Bns$. The $K$-theory of $\Sn$, $\Snp$ and $\Hn$ is still an open problem.

\begin{prop}
 The $K$-groups of the following (full) $C^*$-algebras corresponding to free easy quantum groups are given by:
\begin{itemize}
 \item $K_0(\AOn)=\Z$ and $K_1(\AOn)=\Z$, if $n>2$.
\item $K_0(\ABn)=\Z$ and $K_1(\ABn)=\Z$, if $n>3$.
\item $K_0(\ABnp)=\Z^2$ and $K_1(\ABnp)=\Z^2$, if $n>3$.
\item $K_0(\ABns)=\Z^2$ and $K_1(\ABns)=\Z$, if $n>3$.
\end{itemize}
\end{prop}
\Proof By \cite[th. 9.1]{V11}, we know that $K_0(\AOn)=\Z$ and $K_1(\AOn)=\Z$. In fact, the same holds true for the reduced $C^*$-algebra of $\On$, even in the case of some $Q$-deformed versions of $\On$. This is the main part of the computation using deep theory. The other $K$-groups may be deduced without any effort.

By Proposition \ref{PropRaum}, the $K$-groups of $\ABn$ and $A_{O_{n-1}^+}$ coincide, if $n>3$. 
Since $\ABnp$ is isomorphic to $\ABn\otimes \C^2=\ABn\oplus\ABn$ by Proposition \ref{PropKonstrIsomI}, we immediately deduce the resulting $K$-groups.

For the unital free product of two unital $C^*$-algebras $A$ and $B$, we know the $K$-groups by Cuntz (\cite{Cu82}). This is, if there exist unital homomorphisms $A\to \C$ and $B\to \C$, then $K_0(A*_\C B)$ is the quotient of $K_0(A)\oplus K_0(B)$ by the subgroup generated by $([1_A], -[1_B])$. Secondly, $K_1(A*_\C B)=K_1(A)\oplus K_1(B)$. 
Since $\left(1,-\begin{pmatrix} 1 \\ 1\end{pmatrix}\right)$ is one of the generators of $\Z\oplus \Z^2$, we obtain $K_0(\ABn *_\C C^*(\Z_2))=\Z^2$. Furthermore we compute $K_1(\ABn *_\C C^*(\Z_2))=\Z$.\qqed

Let us end this chapter by the study of exactness of the $C^*$-algebras corresponding to the seven free easy quantum groups. For this, we use the well known result by Simon Wassermann (\cite{Wa76}, \cite{Wa78}, \cite{Wa90}) that the full group $C^*$-algebra $C^*(\F_2)$ of the free group $\F_2$ on two generators is not exact. Furthermore, we use the following statement taken from an article by Paschke and Salinas (\cite[th. 1.1]{PS79}).

\begin{lem}[\cite{PS79}]
 Let $G_1$ and $G_2$ be non-trivial groups and let $G_1$ be of order strictly greater than 2. Then $\F_2$ embeds into the free product $G_1 * G_2$.
\end{lem}
\Proof The generators $x$ and $y$ of $\F_2$ are mapped to $(ac)^2$ and $(bc)^2$, where $a,b\in G_1$ are elements, such that none of the elements $a, b$ and $a^{-1}b$ is the neutral element, and $1\neq c\in G_2$.\qqed

\begin{lem}
 The full group $C^*$-algebra $C^*(\Z_6 * \Z_2)$ of the free product of the groups $\Z_6$ and $\Z_2$ is not exact.
\end{lem}
\Proof By the previous lemma, the free group $\F_2$ on two generators embeds into $\Z_6*\Z_2$. Thus, $C^*(\F_2)$ is a $C^*$-subalgebra of $C^*(\Z_6*\Z_2)$, see for instance \cite[chapter 2, prop. 5.8]{BO08}. Since a $C^*$-subalgebra of an exact $C^*$-algebra is exact again (cf. \cite{KP00}), $C^*(\Z_6 * \Z_2)$ cannot be exact.
\qqed

\begin{thm} 
The $C^*$-algebra $\ASn$ is not exact, if $n\geq 5$.
\end{thm}
\Proof The full group $C^*$-algebra $C^*(\Z_6)$ is isomorphic to $\C^6$ (since $C^*(\Z_6)$ is a six-dimensional, commutative $C^*$-algebra) and thus to $\mathcal C(S_3)$ (since $S_3$ has exactly six elements). This in turn equals $A_{S_3^+}$, since the generators of $A_{S_3^+}$ commute. Secondly, $C^*(\Z_2)$ is isomorphic to $\C^2$ which is isomorphic to $\mathcal C(S_2)=A_{S_2^+}$. We infer that we have a surjection from $A_{S_5^+}$ to $C^*(\Z_6 * \Z_2)=C^*(\Z_6)*_\C C^*(\Z_2)=A_{S_3^+}*_\C A_{S_2^+}$, for the case $n=5$. This can be illustrated by the following matrix:
\[\begin{pmatrix}
A_{S_3^+}  & 0 \\
0   & A_{S_2^+} \\
\end{pmatrix}\]
Since exactness of $C^*$-algebras passes to quotients (see for instance \cite[chapter 10, th. 2.4]{BO08}), $A_{S_5^+}$ is not exact by the above lemma. From the canonical surjections from $A_{S_{n+1}^+}$ to $\ASn$ given by the following matrix, we deduce the result for arbitrary $n\geq 5$.
\[\begin{pmatrix}
\ASn  & 0 \\
0   & 1 \\
\end{pmatrix}\]\qqed

\begin{cor}
 None of the (full) $C^*$-algebras corresponding to the seven free easy quantum groups $\On$, $\Bns$, $\Bnp$, $\Bn$, $\Hn$, $\Snp$, $\Sn$ is exact, if $n\geq 5$. 
\end{cor}
\Proof Again, we use the fact that exactness passes to quotients.
\qqed

\section*{Acknowledgements}

I would like to thank Roland Speicher for many helpful discussions.
I am indebted to Octavio Arizmendi  for his computations of the moments  of Proposition \ref{RemOctavio}. Furthermore, I would like to thank Michael Brannan for the discussions on exactness of easy quantum groups. Also, I am indebted to Christian Voigt for some remarks on the computation of the $K$-groups in this article. Finally, I would like to thank Teo Banica for his help to improve this article.

\end{document}